\documentclass[reqno,11pt]{amsart}

\usepackage[export]{adjustbox}
\usepackage[foot]{amsaddr}
\usepackage{array}
\usepackage{bbold}
\usepackage[font=small,skip=0pt]{caption}
\usepackage{charter}
\usepackage{enumitem}
\usepackage{float}
\usepackage[margin=1in]{geometry}
\usepackage[colorlinks=true,linktoc=all,linkcolor=blue,citecolor=blue]{hyperref}
\usepackage{multicol}
\usepackage[most]{tcolorbox}
\usepackage{tikz}

\newcolumntype{x}[1]{>{\centering\arraybackslash\hspace{0pt}}p{#1}}

\newtcolorbox[auto counter]{picbox}{%
	enhanced,
	width=\linewidth,
	boxrule=.5pt,
	colback=white,
	colupper=black,
	sharp corners,
	lower separated=true,
}

\theoremstyle{plain}
\newtheorem{theorem}{Theorem}[section]
\numberwithin{equation}{section}

\usetikzlibrary{arrows}
\tikzstyle{state} = [rectangle, draw, node distance=1.75cm, text width=2em, text centered, minimum height=2em, thick]

\title[Sociological SIR Models]{Analysis of SIR Epidemic Models with Sociological Phenomenon}
\author{Robert F.~Allen\textsuperscript{1}, Katherine C.~Heller\textsuperscript{2}, and Matthew A.~Pons\textsuperscript{2}}
\address{\textsuperscript{1}Department of Mathematics \& Statistics, University of Wisconsin-La Crosse}
\address{\textsuperscript{2}Department of Mathematics \& Actuarial Science, North Central College}
\email{rallen@@uwlax.edu, kheller@noctrl.edu, mapons@nocrtl.edu}


\begin{document}

\begin{abstract}
We propose two SIR models which incorporate sociological behavior of groups of individuals.  It is these differences in behaviors which impose different infection rates on the individual susceptible populations, rather than biological differences.  We compute the basic reproduction number for each model, as well as analyze the sensitivity of $R_0$ to changes in sociological parameter values.  
\end{abstract}

\maketitle

\section{Introduction}
The use of SIR models in the mathematical study of communicable diseases began in 1927 with the work of Kermack and McKendrick \cite{KermackMcKendrick:1927}.  Since then, SIR models, and their variants, have been used to model different population structures for different diseases.  For an excellent reference on Mathematical Epidemiology see \cite{Brauer:2017}. Traditionally, SIR models have been employed to study biological aspects of disease propagation.  This includes multiple populations with varying immunities, either natural or vaccine induced, see for example \cite{LiuJiangHyatAhmad:2018, MengChen:2008, SahuDhar:2012}.

In recent years, sociological interventions have been employed to reduce infection rates of global diseases, such as Severe Acute Respiratory Syndrome (SARS), H1N1, and most recently COVID-19. Many such sociological interventions include self-imposed quarantine.  SIR models with such quarantine interventions have been devised and analyzed.  The interested reader is directed to \cite{AlamKabirTanimoto:2020, AllenJensWendt:2014, HethcoteZhienShengbing:2002} as examples.   
In this paper, we propose two SIR models which include sociological differences among 2 susceptible populations.  These sociological differences impose different rates of infection, rather than biological differences.  Examples of such sociological differences might include
\begin{itemize}
\item front-line workers versus general population during a global pandemic
\item pro-maskers versus anti-maskers
\item individuals that adhere to a stay-at-home suggestion and those that do not. 
\end{itemize} While many of the situations that come to mind can be linked to COVID-19, the intent of this paper is not to delve into the political use of such a model.  As the models include sociological behavior, data does not necessarily exist to validate the models against or estimate parameters.  As what will be discussed in Sections \ref{Section:Conclusion} and \ref{Section:Future}, the power of these models will be in the predicting how different sociological interventions can effect the spread of a disease.  


\section{Simplified Sociological Model}\label{Section:BasicModel}

In this section, we devise a mathematical model, referred to as \eqref{System:BasicModel}, of two sub-populations $S_1$ and $S_2$ of healthy individuals.  Membership in each sub-population is determined by sociological and not biological reasons.  The groups $S_1$ and $S_2$ are biologically identical, in that individuals in either group are equally likely to contract the disease and recover identically.  However, the behavior of each group determines differences in infection rates.  The behavior of $S_1$ is such that the incident of infection is greater than that of $S_2$.  These rates are denoted by $\beta_1$ and $\beta_2$.  To model different types of group behavior, we assume $\beta_1 > \beta_2$, so the two groups are distinguishable.  
In this model, we have two infected sub-populations, those infected that exhibit symptoms and those who are asymptomatic.  In this model, any infected person is infectious.  The symptomatic group is denoted by $I_s$ and the asymptomatic group by $I_a$.  Any infectious individual can contribute to the spread of the disease, and thus the infectious individuals are denoted by $I = I_a+I_s$.  We assume a mass-action transmission of the disease between the groups $S_1$, $S_2$, and $I$.  Since the healthy populations $S_1$ and $S_2$ are biologically identical, the proportion of infections that result in symptomatic individuals is denoted by $\lambda$, where the proportion of infections resulting in an asymptomatic individual is $1-\lambda$.  

The model we propose in this section allows for recovery, which we denote by the group $R$.  We assume the rate of recovery is the same for individuals in either $I_a$ or $I_s$, which we denote by $\kappa$.  Likewise, asymptomatic individuals develop symptoms at a rate of $\gamma$.   Finally, we assume our population to be closed with a total population $N = S_1 + S_2 + I_a + I_s + R$. The component diagram for \eqref{System:BasicModel} is provided in Figure \ref{Figure:AsymptomatircTwoSusceptibleModel}, and the system of ordinary differential equations is given below.

\begin{figure}[H]
	\begin{center}
		\begin{tikzpicture}[->,>=stealth',every node/.style={font=\footnotesize}]
			\coordinate (center);
			\node[state] (S1) [left of=center, node distance=3cm] {$S_1$};
			\node[state] (Is) [above of=center] {$I_s$};
			\node[state] (Ia) [below of=center] {$I_a$};
			\node[state] (R)  [right of=center, node distance=1cm] {$R$};
			\node[state] (S2) [right of=R, node distance=2cm] {$S_2$};
			
			\path 
			(S1) edge [thick] node [below, rotate=-30] {$(1-\lambda)\textcolor{red}{\beta_1}S_1I$} (Ia)
			(S1) edge [thick] node [above, rotate=30] {$\lambda\textcolor{red}{\beta_1}S_1I$} (Is)
			(S2) edge [thick] node [below, rotate=30] {$(1-\lambda)\textcolor{red}{\beta_2}S_2I$} (Ia)
			(S2) edge [thick] node [above, rotate=-30] {$\lambda\textcolor{red}{\beta_2}S_2I$} (Is)
			(Ia) edge [thick, dashed, bend left] node [left] {$\gamma I_a$} (Is)
			(Is) edge [thick, dashed] node [right] {$\kappa I_s$} (R)
			(Ia) edge [thick, dashed] node [right] {$\kappa I_a$} (R);
		\end{tikzpicture}
	\end{center}
	\caption{Compartment diagram of simplified two susceptible population model.}\label{Figure:AsymptomatircTwoSusceptibleModel}
\end{figure}
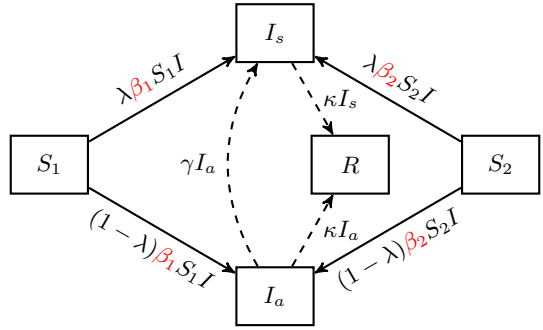

{\everymath={\displaystyle}
	\begin{equation}\label{System:BasicModel}\tag{MA}
		\begin{array}{c@{}l}
			\frac{dS_1}{dt} &{}= -\beta_1 \frac{S_1}{N} (I_a+I_s)\\[1em]
			\frac{dS_2}{dt} &{}= -\beta_2 \frac{S_2}{N} (I_a+I_s)\\[1em]
			\frac{dI_s}{dt} &{}= \lambda\left[\beta_1\frac{S_1}{N} + \beta_2\frac{S_2}{N}\right](I_a+I_s) + \gamma I_a - \kappa I_s\\[1em]
			\frac{dI_a}{dt} &{}= (1-\lambda)\left[\beta_1\frac{S_1}{N} + \beta_2\frac{S_2}{N}\right](I_a+I_s) - (\gamma + \kappa)I_a\\[1em]
			\frac{dR}{dt} &{}= \kappa(I_a + I_s)
		\end{array}
\end{equation}}

Appendix B provides several numerical simulations of model \eqref{System:BasicModel}, with an eye toward the effect of changes in the system parameters $\beta_1, \beta_2,$ and $\rho$.  Specifically, the reader is directed to Figures \ref{Figure:BasicModelSimulations}, \ref{Figure:BasicModelEffectOfB1}, \ref{Figure:BasicModelEffectOfB2}, and \ref{Figure:BasicModelChangeinRho}.

\subsection{Stability Analysis}
We first compute $R_0$ using the Next Generation Matrix, as outlined in \cite{DiekmannHeesterbeekRoberts:2010}.  For the disease-free equilibrium, there exists a $0 < \rho < 1$ such that $S_1^* = \rho N$ and thus $S_2^* = (1-\rho)N$.  The linearized infection subsystem, about the disease-free equilibrium $(\rho N, (1-\rho)N, 0, 0, 0)$, where \[B_\rho = \beta_1\rho + \beta_2(1-\rho)\] is given by 

{\everymath={\displaystyle}
	\[
	\begin{array}{c@{}l}
		\frac{dI_s}{dt} &{}= \lambda B_\rho(I_a+I_s) + \gamma I_a - \kappa I_s \\[1em]
		\frac{dI_a}{dt} &{}= (1-\lambda) B_\rho(I_a+I_s) - (\gamma + \kappa)I_a
	\end{array}
	\]
} 
	
\noindent and the Jacobian of the linearized infection subsystem by
\[
J =
\left(\begin{matrix}
	\lambda B_\rho - \kappa & \lambda B_\rho + \gamma\\[.5em]
	(1-\lambda) B_\rho & (1-\lambda) B_\rho - (\gamma+\kappa)
\end{matrix}\right).
\] We decompose the Jacobian matrix as $J = T + \Sigma$, where $T$ is the new infection matrix and $\Sigma$ is the transition matrix.  Thus 
\[T = \left(\begin{matrix}
	\lambda B_\rho & \lambda B_\rho\\[.5em]
	(1-\lambda) B_\rho & (1-\lambda) B_\rho
\end{matrix}\right)\]  and \[\Sigma = 
\left(\begin{matrix}
	-\kappa & \gamma\\[.5em]
	0 & -(\gamma+\kappa)
\end{matrix}\right).\] The Next Generation Matrix is calculated as
\[K = -T\Sigma^{-1\\} = \frac{1}{\kappa}B_\rho
\left(\begin{matrix}
	\lambda & \lambda \\
	1-\lambda & 1-\lambda
\end{matrix}\right).\]

The eigenvalues of $K$ are 0 and
\[\ell = \frac{1}{\kappa}\left[\beta_1\rho + \beta_2(1-\rho)\right].\]  Since all system parameters are positive and $\rho \in (0,1)$, it follows that $\ell > 0$, and thus the dominant eigenvalue of $K$.  

\begin{theorem}\label{Theorem:BaseModel:R0}
	For system \eqref{System:BasicModel}, the basic reproduction number is
	\[R_0 = \frac{1}{\kappa}\left[\rho\beta_1 + (1-\rho)\beta_2\right].\] Moreover, a disease-free equilibrium of the form $(\rho N, (1-\rho)N, 0, 0, 0)$ for $\rho \in (0,1)$ is asymptotically stable if and only if \begin{equation}\label{Equation:BaseModel:R0ConvexityCondition}
		\rho\beta_1 + (1-\rho)\beta_2 < \kappa.
	\end{equation}
\end{theorem}
	
From Theorem \ref{Theorem:BaseModel:R0}, we see that $R_0$ is a function of $\kappa, \rho, \beta_1,$ and $\beta_2$.  Since $\kappa$ is assumed to be constant, and dependent only on the disease, we will analyze $R_0$ as a function of $\rho, \beta_1,$ and $\beta_2$.  First, we will determine for a fixed value of $\rho$, the values of $\beta_1$ and $\beta_2$ for which $R_0 < 1$.  Such a pair $(\beta_1,\beta_2)$ is called $\rho$-feasible. 

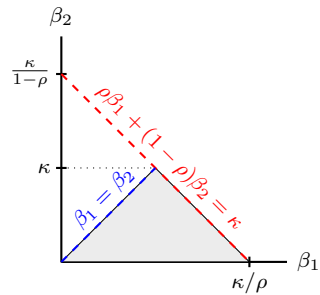
\begin{figure}[H]
	\begin{center}
		\begin{tikzpicture}
			\draw [thick] (0,0)-- (0,3);
			\draw [thick] (0,0) -- (3,0);
			\draw [thick] (2.5,.07) -- (2.5,-.07) node [below] at (2.5,0) {\scriptsize$\kappa/\rho$};
			\draw [thick] (-.07,2.5) -- (.07,2.5) node [left] at (0,2.5) {\scriptsize$\frac{\kappa}{1-\rho}$};
			\draw [fill=gray!15] (0,0) -- (1.25,1.25) -- (2.5,0) -- cycle;
			\draw [thick, dashed, red] (0,2.5) -- (2.5,0);
			\draw [thick, dashed, blue] (0,0) -- (1.25,1.25);
			\draw [dotted] (0,1.25) -- (1.25,1.25);
			\draw [thick] (.07,1.25) -- (-.07,1.25) node [left] at (0,1.25) {\scriptsize$\kappa$};
			
			\draw node [above] at (0,3) {\scriptsize$\beta_2$};
			\draw node [right] at (3,0) {\scriptsize$\beta_1$};
			\draw node [above, rotate=45, blue] at (.75,.675) {\scriptsize$\beta_1=\beta_2$};
			\draw node [above, rotate=-45, red] at (1.25,1.2) {\scriptsize$\rho\beta_1 + (1-\rho)\beta_2 = \kappa$};
		\end{tikzpicture}
		\caption{Possible $\rho$-feasible set.}\label{Figure:BaseModel:PossibleFeasibleSet}
	\end{center}
\end{figure} 

First note that it is assumed in the model that $\beta_1 > \beta_2$.  Also, $(\beta_1,\beta_2)$ is $\rho$-feasible if and only if $\rho\beta_1 + (1-\rho)\beta_2 < \kappa$ from inequality \eqref{Equation:BaseModel:R0ConvexityCondition}.  Thus, the set of $\rho$-feasible points must satisfy the system of linear inequalities
\[\begin{cases}
	\begin{aligned}
		\beta_1 &> \beta_2,\\
		\kappa &> \rho\beta_1 + (1-\rho)\beta_2
	\end{aligned}
\end{cases}
\] 
and thus must be a subset of the shaded region in Figure \ref{Figure:BaseModel:PossibleFeasibleSet}.  More specifically, the $\rho$-feasible set is a subset of the convex hull of the points $\{(0,0), (\kappa,\kappa), (\kappa/\rho,0)\}$.

The following is an immediate consequence.

\begin{theorem}\label{Theorem:BasicModel|StabilityCriterion1} The disease-free equilibrium $(\rho N, (1-\rho)N, 0, 0, 0)$ of System \eqref{System:BasicModel} is asymptotically unstable (i.e.~$R_0 \geq 1$) if $\beta_2 \geq \kappa$.
\end{theorem}
If we impose the constraint that $\beta_2 < \kappa$, then the set of $\rho$-feasible points for System \eqref{System:BasicModel} must satisfy the linear inequalities
\[\begin{cases}
	\begin{aligned}
		\beta_1 &> \beta_2,\\
		\kappa &> \beta_2,\\
		\kappa &> \rho\beta_1 + (1-\rho)\beta_2.
	\end{aligned}
\end{cases}
\] 
The following linear stability of disease-free equilibrium are established, each determining the $\rho$-feasible points given two of the three quantities $\rho,\beta_1,\beta_2$ fixed.

\begin{theorem}\label{Theorem:BasicModel|StabilityCriterion2} For fixed $\beta_1 \in (0,1)$ and $0 < \beta_2 < \min\{\beta_1, \kappa\},$ the disease-free equilibrium $(\rho N, (1-\rho)N, 0, 0, 0)$ of System \eqref{System:BasicModel} is asymptotically stable if and only if $\rho \in (0,\mathrm{P})$ where 
	\[\mathrm{P} = \frac{\kappa-\beta_2}{\beta_1-\beta_2}.\]     
\end{theorem}

\begin{theorem}\label{Theorem:BasicModel|StabilityCriterion3} For fixed $\rho,\beta_1 \in (0,1)$, the disease-free equilibrium $(\rho N, (1-\rho)N,0,0,0)$ of System \eqref{System:BasicModel} is asymptotically stable if and only if $\beta_2 \in (0,\mathrm{B}_2)$ where \[\mathrm{B}_2 = \begin{cases}\beta_1 & \text{ \normalfont{if} $\beta_1 \leq \kappa$}\\\displaystyle\frac{\kappa-\rho\beta_1}{1-\rho} & \text{ \normalfont{if} $\kappa < \beta_1 \leq \min\left\{1,\frac{\kappa}{\rho}\right\}$}.\end{cases}\]
\end{theorem}

\begin{theorem}\label{Theorem:BasicModel|StabilityCriterion4} For fixed $\rho \in (0,1)$ and $\beta_2 \in (0,\kappa)$, the disease-free equilibrium $(\rho N,(1-\rho)N,0,0,0)$ of System \eqref{System:BasicModel} is asymptotically stable if and only if $\beta_1 \in (\beta_2,\mathrm{B}_1]$ where \[\mathrm{B}_1 = \min\left\{1,\frac{\kappa-(1-\rho)\beta_2}{\rho}\right\}.\]
\end{theorem}

As discussed in \cite[Section 9.2]{BrauerCastillo-Chavez:2012}, the number of contacts made by an individual to infect a member of the $S_1$ population is $\beta_1 N$.  Thus, it is natural to take $\beta_1 \in (0,1]$ and $\beta_2 \in (0,1)$, since $\beta_2 < \beta_1$.  Then the set of $\rho$-feasible points will be the intersection of the shaded region in Figure \ref{Figure:BaseModel:PossibleFeasibleSet} with the set $(0,1]\times (0,1)$.  

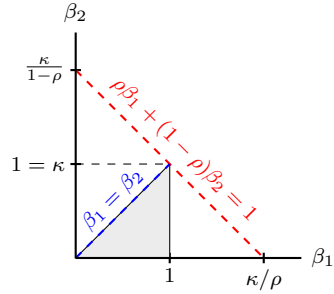
\begin{figure}[H]
	\begin{center}
		\begin{tikzpicture}
			\draw [thick] (1.25,.07) -- (1.25,-.07) node [below] at (1.25,0) {\scriptsize$1$};
			\draw [thick] (-.07,1.25) -- (.07,1.25) node [left] at (0,1.25) {\scriptsize$1=\kappa$};
			\draw [dashed] (0,1.25) -- (1.25,1.25);
			\draw [fill=gray!15] (0,0) -- (1.25,1.25) -- (1.25,0) -- cycle;
			\draw [thick] (2.5,.07) -- (2.5,-.07) node [below] at (2.5,0) {\scriptsize$\kappa/\rho$};
			\draw [thick] (-.07,2.5) -- (.07,2.5) node [left] at (0,2.5) {\scriptsize$\frac{\kappa}{1-\rho}$};
			\draw [thick, dashed, red] (0,2.5) -- (2.5,0);
			\draw [thick, dashed, blue] (0,0) -- (1.25,1.25);
			
			\draw node [above] at (0,3) {\scriptsize$\beta_2$};
			\draw node [right] at (3,0) {\scriptsize$\beta_1$};
			\draw node [blue, above, rotate=45] at (.675,.6) {\scriptsize$\beta_1=\beta_2$};
			\draw node [red, above, rotate=-45] at (1.25,1.25) {\scriptsize$\rho\beta_1 + (1-\rho)\beta_2 = 1$};
			\draw [thick] (0,0)-- (0,3);
			\draw [thick] (0,0) -- (3,0);
		\end{tikzpicture}
		\caption{Feasible sets of System \eqref{System:BasicModel} for $\kappa = 1$.}\label{Figure:BaseModel:ActualFeasibleSetsKLarge}
	\end{center}
\end{figure}

At any time, $\kappa(I_a+I_s)$ is the number of new recovered individuals.  Thus, it must be the case that $\kappa \in (0,1]$, or else the number of new recovered individuals will exceed the number of infected individuals which is biologically impossible.  If $\kappa = 1$, then the line $\rho\beta_1 + (1-\rho)\beta_2 = \kappa$ intersects with the set $(0,1]\times (0,1)$ at the point $(1,1)$.  Thus, for every $\rho \in (0,1)$, every $(\beta_1,\beta_2) \in (0,1]\times (0,1)$ which satisfies $\beta_2 < \beta_1$ is $\rho$-feasible, as depicted in Figure \ref{Figure:BaseModel:ActualFeasibleSetsKLarge}.

We now consider the case when $\kappa < 1$ to determine the $\rho$-feasible points $(\beta_1,\beta_2)$.  For this case, the line $\rho\beta_1 + (1-\rho)\beta_2 = \kappa$ will intersect non-trivially with the set $(0,1]\times (0,1)$.  However, the exact shape of the resulting $\rho$-feasibility set will be determined by the location of the intercepts $\frac{\kappa}{\rho}$ and $\frac{\kappa}{1-\rho}$ with respect to the unit square $[0,1]\times [0,1]$. 

\begin{figure}[h]
	\begin{center}
		\begin{tabular}{c c c}
			\begin{tikzpicture}
				\fill [gray!20] (0,0) -- (1.25, 1.25) -- (2.5,.5) -- (2.5,0) -- (0,0) -- cycle;
				\draw [dashed, thick] (0,0) -- (0,2.5);
				\draw [thick] (0,2.5) -- (0,2.75);
				\draw node [above] at (0,2.75) {\scriptsize$\beta_2$};
				\draw [dashed, thick] (0,0) -- (2.5,0);
				\draw [thick] (2.5,0) -- (2.75,0);
				\draw node [right] at (2.75,0) {\scriptsize$\beta_1$};
				\draw [thick] (2.5,.07) -- (2.5,-.07) node [below] at (2.5,0) {\scriptsize$1$};
				\draw [thick] (.07,2.5) -- (-.07,2.5) node [left] at (0,2.5) {\scriptsize$1$};
				\draw [dashed, thick] (0,2.5) -- (2.5,2.5);
				\draw [thick] (2.5,0) -- (2.5,2.5);
				\draw [dashed, thick, blue] (0,0) -- (2.5,2.5);
				\draw [dashed, thick, red] (0,2) -- (2.5,.5);
				\draw [thick] (-.07,2) -- (.07,2) node [left] at (-.1,2) {\scriptsize$\frac{\kappa}{1-\rho}$};
				\draw [thick] (2.43,.5) -- (2.57,.5)  node [right] at (2.5,.7) {\scriptsize$\frac{\kappa-\rho}{1-\rho}$};
				\draw [dotted] (0,1.25) -- (1.25,1.25);
				\draw [thick] (-.07,1.25) -- (.07,1.25) node [left] at (0,1.25) {\scriptsize$\kappa$};
				\draw node [below] at (2,0) {\textcolor{white}{\scriptsize$\kappa/\rho$}};
				\draw node [above] at (1,2.5) {\textcolor{white}{\scriptsize$\frac{\kappa-(1-\rho)}{\rho}$}};
				
			\end{tikzpicture}
			&
			\begin{tikzpicture}
				\fill [gray!20] (0,0) -- (1.25, 1.25) -- (2.5,0) -- (0,0) -- cycle;
				\draw [dashed, thick] (0,0) -- (0,2.5);
				\draw [thick] (0,2.5) -- (0,2.75);
				\draw node [above] at (0,2.75) {\scriptsize$\beta_2$};
				\draw [dashed, thick] (0,0) -- (2.5,0);
				\draw [thick] (2.5,0) -- (2.75,0);
				\draw node [right] at (2.75,0) {\scriptsize$\beta_1$};
				\draw [thick] (2.5,.07) -- (2.5,-.07) node [below] at (2.5,0) {\scriptsize$1$};
				\draw [thick] (.07,2.5) -- (-.07,2.5) node [left] at (0,2.5) {\scriptsize$1$};
				\draw [dashed, thick] (0,2.5) -- (2.5,2.5);
				\draw [thick] (2.5,0) -- (2.5,2.5);
				\draw [dashed, thick, blue] (0,0) -- (2.5,2.5);
				\draw [dashed, thick, red] (0,2.5) -- (2.5,0);
				\draw [dotted] (0,1.25) -- (1.25,1.25);
				\draw [thick] (-.07,1.25) -- (.07,1.25) node [left] at (0,1.25) {\scriptsize$\kappa$};
				\draw node [below] at (2,0) {\textcolor{white}{\scriptsize$\kappa/\rho$}};
				\draw node [right] at (2.5,.5) {\textcolor{white}{\scriptsize$\frac{\kappa-\rho}{1-\rho}$}};
				\draw node [above] at (1,2.5) {\textcolor{white}{\scriptsize$\frac{\kappa-(1-\rho)}{\rho}$}};
				
			\end{tikzpicture}
			& 
			\begin{tikzpicture}
				\fill [gray!20] (0,0) -- (1.25, 1.25) -- (2,0) -- (0,0) -- cycle;
				\draw [dashed, thick] (0,0) -- (0,2.5);
				\draw [thick] (0,2.5) -- (0,2.75);
				\draw node [above] at (0,2.75) {\scriptsize$\beta_2$};
				\draw [dashed, thick] (0,0) -- (2.5,0);
				\draw [thick] (2.5,0) -- (2.75,0);
				\draw node [right] at (2.75,0) {\scriptsize$\beta_1$};
				\draw [thick] (2.5,.07) -- (2.5,-.07) node [below] at (2.5,0) {\scriptsize$1$};
				\draw [thick] (.07,2.5) -- (-.07,2.5) node [left] at (0,2.5) {\scriptsize$1$};
				\draw [dashed, thick] (0,2.5) -- (2.5,2.5);
				\draw [thick] (2.5,0) -- (2.5,2.5);
				\draw [dashed, thick, blue] (0,0) -- (2.5,2.5);
				\draw [dashed, thick, red] (.5,2.5) -- (2,0);
				\draw [thick] (.5,2.57) -- (.5,2.43) node [above] at (1,2.5)
				{\scriptsize$\frac{\kappa-(1-\rho)}{\rho}$};
				\draw [thick] (2,-.07) -- (2,.07) node [below] at (2,0) {\scriptsize $\kappa/\rho$};
				\draw [dotted] (0,1.25) -- (1.25,1.25);
				\draw [thick] (-.07,1.25) -- (.07,1.25) node [left] at (0,1.25) {\scriptsize$\kappa$};
				\draw node [right] at (2.5,.5) {\textcolor{white}{\scriptsize$\frac{\kappa-\rho}{1-\rho}$}};
			\end{tikzpicture}\\[-.65em]
			\hspace{-.2in} {\small$\rho < \kappa$} & \hspace{-.2in} {\small$\rho = \kappa$} & \hspace{-.2in} {\small$\rho > \kappa$}\\[.5em]
			\hspace{-.2in} Type 1 & \hspace{-.2in} Type 0 & \hspace{-.2in} Type $-1$
		\end{tabular}
		\caption{Feasible sets of System \eqref{System:BasicModel} for $\kappa = \frac{1}{2}$.}\label{Figure:BaseModel:ActualFeasibleSetsK1/2}
	\end{center}
\end{figure}
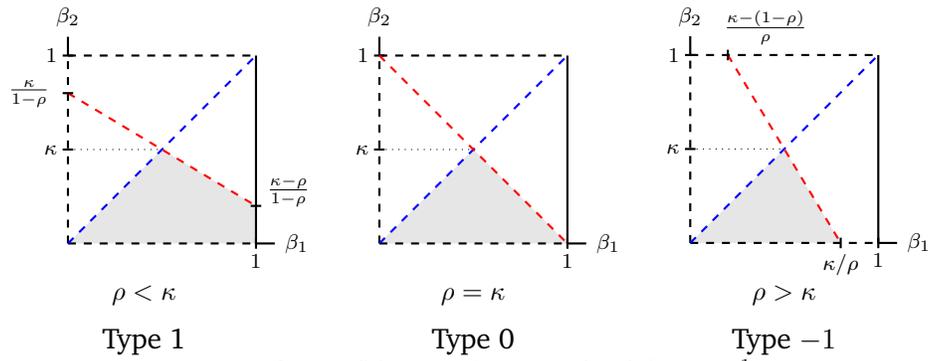

First, we consider the case when $\kappa = \frac12$. The set of $\rho$-feasible points is depicted in Figure \ref{Figure:BaseModel:ActualFeasibleSetsK1/2} and the $\rho$-feasible sets in the cases of $0 < \kappa < \frac12$ and $\frac12 < \kappa < 1$ are shown in Figure \ref{Figure:ActualFeasibleSets}. The notable differences between the feasible set types is the existence of $(\beta_1,\beta_2)$ pairs which guarantee an asymptotically stable disease-free equilibrium. From these figures, we classify $\rho$-feasible sets as follows:
\begin{description}
	\item[Type 1:] The $\rho$-feasible set is the convex hull of the four distinct points $\{(0,0), (\kappa,\kappa), (1,(\kappa-\rho)/(1-\rho)), (1,0)\}$.
	
	For every $\beta_1 \in (0,1]$, there exists a $0 < \beta_2 < \beta_1$ for which $(\beta_1,\beta_2)$ is $\rho$-feasible.  Specifically, there exist $\rho$-feasible pairs of the form $(1,\beta_2)$.  Thus, for every $\rho \in (0,1)$ and every $\beta_1 \in (0,1]$, there exists a $\beta_2$ for which the disease-free equilibrium is asymptotically stable.
	
	\item[Type 0:] The $\rho$-feasible set is the convex hull of the three distinct points $\{(0,0), (\kappa,\kappa), (1,0)\}$.
	
	The only $(\beta_1,\beta_2)$ pairs that are not $\rho$-feasible are of the form $(1,\beta_2)$.  So for every $\rho \in (0,1)$ and $\beta_1$ in $(0,1)$, there exists $\beta_2$ for which the disease-free equilibrium is linearly stable.
	
	\item[Type -1:] The $\rho$-feasible set is the convex hull of the three distinct points $\{(0,0), (\kappa,\kappa), (\kappa/\rho,0)\}$. 
	
	There are no $(\beta_1,\beta_2)$ pairs that are $\rho$-feasible if $\frac{\kappa}{\rho} < \beta_1 \leq 1$.  So for every $\rho \in (0,1)$ and every $\beta_1 \geq \frac{\kappa}{\rho}$ there are no values of $\beta_2$ for which the disease-free equilibrium is asymptotically stable.
\end{description}
In fact, we see that for every value of $\rho \in (0,1)$, the convex hull of the points $\{(0,0), (\kappa,\kappa), (\kappa,0)\}$ is contained in the feasible set.  

\begin{theorem} For System \eqref{System:BasicModel}, there exist $\beta_1$ and $\beta_2$ for which the disease-free equilibrium $(\rho N, (1-\rho)N, 0, 0, 0)$ is asymptotically stable for every $\rho \in (0,1)$. Specifically, if $\beta_1 < \kappa$, then for all $\beta_2 \in (0,\beta_1)$, the point $(\beta_1,\beta_2)$ is $\rho$-feasible for every $\rho \in (0,1)$.
\end{theorem}

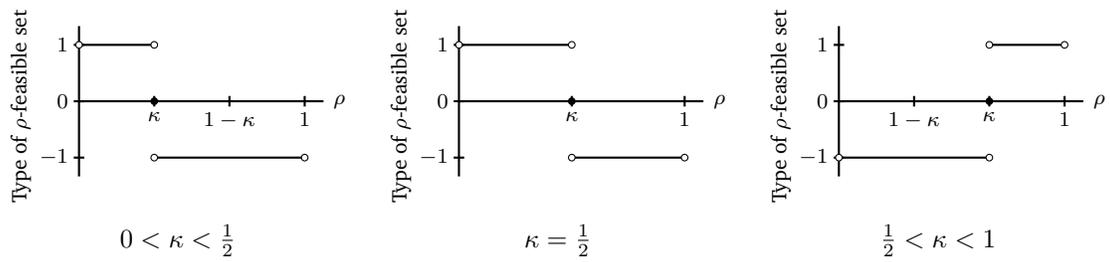
\begin{figure}[h]
	\begin{center}
		\begin{tabular}{c c c}
			\begin{tikzpicture}
				\draw [thick] (0,0) -- (3.25,0);
				\draw [thick] (0,-1) -- (0,1);
				\draw node [right] at (3.25,0) {\scriptsize$\rho$};
				\draw node [left, rotate=90] at (-.77,1.35) {\scriptsize Type of $\rho$-feasible set};
				\draw [thick] (1,.07) -- (1,-.07) node [below] at (1,0) {\scriptsize$\kappa$};
				\draw [thick] (2,.07) -- (2,-.07) node [below] at (2,0) {\scriptsize$1-\kappa$};
				\draw [thick] (3,.07) -- (3,-.07) node [below] at (3,0) {\scriptsize$1$};
				\draw [thick] (.07,.75) -- (-.07,.75) node [left] at (0,.75) {\scriptsize$1$};
				\draw [thick] (.07,-.75) -- (-.07,-.75) node [left] at (0,-.75) {\scriptsize$-1$};
				\draw [thick] (.07,0) -- (-.07,0) node [left] at (0,0) {\scriptsize$0$};
				\draw [thick] (0,.75) -- (1,.75);
				\draw [thick] (1,-.75) -- (3,-.75);
				\filldraw [fill=white] (0,.75) circle (1.25pt);
				\filldraw [fill=white] (1,.75) circle (1.25pt);
				\filldraw [fill=white] (1,-.75) circle (1.25pt);
				\filldraw [fill=white] (3,-.75) circle (1.25pt);
				\filldraw (1,0) circle (1.25pt);
			\end{tikzpicture}
			&
			\begin{tikzpicture}
				\draw [thick] (0,0) -- (3.25,0);
				\draw [thick] (0,-1) -- (0,1);
				\draw node [right] at (3.25,0) {\scriptsize$\rho$};
				\draw node [left, rotate=90] at (-.77,1.35) {\scriptsize Type of $\rho$-feasible set};
				\draw [thick] (1.5,.07) -- (1.5,-.07) node [below] at (1.5,0) {\scriptsize$\kappa$};
				\draw [thick] (3,.07) -- (3,-.07) node [below] at (3,0) {\scriptsize$1$};
				\draw [thick] (.07,0) -- (-.07,0) node [left] at (0,0) {\scriptsize$0$};
				\draw [thick] (.07,.75) -- (-.07,.75) node [left] at (0,.75) {\scriptsize$1$};
				\draw [thick] (.07,-.75) -- (-.07,-.75) node [left] at (0,-.75) {\scriptsize$-1$};
				\draw [thick] (0,.75) -- (1.5,.75);
				\draw [thick] (1.5,-.75) -- (3,-.75);
				\filldraw [fill=white] (0,.75) circle (1.25pt);
				\filldraw [fill=white] (1.5,.75) circle (1.25pt);
				\filldraw [fill=white] (1.5,-.75) circle (1.25pt);
				\filldraw [fill=white] (3,-.75) circle (1.25pt);
				\filldraw (1.5,0) circle (1.25pt);
			\end{tikzpicture}
			&
			\begin{tikzpicture}
				\draw [thick] (0,0) -- (3.25,0);
				\draw [thick] (0,-1) -- (0,1);
				\draw node [right] at (3.25,0) {\scriptsize$\rho$};
				\draw node [left, rotate=90] at (-.77,1.35) {\scriptsize Type of $\rho$-feasible set};
				\draw [thick] (1,.07) -- (1,-.07) node [below] at (1,0) {\scriptsize$1-\kappa$};
				\draw [thick] (2,.07) -- (2,-.07) node [below] at (2,0) {\scriptsize$\kappa$};
				\draw [thick] (3,.07) -- (3,-.07) node [below] at (3,0) {\scriptsize$1$};
				\draw [thick] (.07,.75) -- (-.07,.75) node [left] at (0,.75) {\scriptsize$1$};
				\draw [thick] (.07,-.75) -- (-.07,-.75) node [left] at (0,-.75) {\scriptsize$-1$};
				\draw [thick] (.07,0) -- (-.07,0) node [left] at (0,0) {\scriptsize$0$};
				\draw [thick] (0,-.75) -- (2,-.75);
				\draw [thick] (2,.75) -- (3,.75);
				\filldraw [fill=white] (0,-.75) circle (1.25pt);
				\filldraw [fill=white] (2,-.75) circle (1.25pt);
				\filldraw [fill=white] (2,.75) circle (1.25pt);
				\filldraw [fill=white] (3,.75) circle (1.25pt);
				\filldraw (2,0) circle (1.25pt);
			\end{tikzpicture}\\
			{\small$0<\kappa<\frac12$} & {\small$\kappa=\frac12$} & {\small$\frac12 < \kappa < 1$}
		\end{tabular}
	\end{center}
	\caption{Bifurcation diagrams of System \eqref{System:BasicModel}  as a function of $\rho$.}\label{Figure:BifurcationTwoSusceptibleModel}
\end{figure}

As $\rho$ varies from 0 to 1, the System \eqref{System:BasicModel} undergoes a bifurcation in the type of $\rho$-feasible set which guarantees asymptotic stability of disease-free equilibria, as depicted in Figure \ref{Figure:BifurcationTwoSusceptibleModel}.  One can also view Figure \ref{Figure:BifurcationTwoSusceptibleModel} as a bifurcation in the feasible type profile of System \eqref{System:BasicModel} as the parameter $\kappa$ varies from 0 to 1.  

\subsection{Sensitivity Analysis}\label{Subsection:BasicModelSentitivityAnalysis}
In \cite{ChitnisHymanCushing:2008}, the normalized forward sensitivity index of a variable $u$ that depends differentiably on a parameter $p$ is defined as \[\Upsilon_p^u = \left(\frac{p}{u}\right)\frac{\partial u}{\partial p}.\]  This sensitivity index is a measure of the relative change in a state variable with respect to a parameter.  In this section, we compute the sensitivity indices of $R_0$ with respect to the three sociological parameters $\rho, \beta_1$, and $\beta_2$.  

From Theorem \ref{Theorem:BaseModel:R0}, we have the explicit expression for $R_0$ as \[R_0 = \frac{1}{\kappa}\left[\rho\beta_1 + (1-\rho)\beta_2\right]\] or \[R_0\kappa = \rho\beta_1 + (1-\rho)\beta_2.\]  Observe all sensitivity indices of $R_0$ with respect to the sociological parameters $\rho, \beta_1$, and $\beta_2$ have the same structure of $1-\frac{f(\rho,\beta_1,\beta_2)}{R_0\kappa}$. 

\begin{equation}\label{Equation:MASensitivityIndices}
\begin{aligned}
\Upsilon_\rho^{R_0} &= \left(\frac{\rho}{R_0}\right)\frac{\partial R_0}{\partial\rho} = 1-\frac{\beta_2}{R_0\kappa}\\
\Upsilon_{\beta_1}^{R_0} &= \left(\frac{\beta_1}{R_0}\right)\frac{\partial R_0}{\partial \beta_1} = 1-\frac{(1-\rho)\beta_2}{R_0\kappa}\\
\Upsilon_{\beta_2}^{R_0} &= \left(\frac{\beta_2}{R_0}\right)\frac{\partial R_0}{\partial \beta_2} = 1-\frac{\rho\beta_1}{R_0\kappa}
\end{aligned}
\end{equation}

Immediately, we see that each sensitivity index is positive, less than 1, and $\Upsilon_{\beta_1}^{R_0} + \Upsilon_{\beta_2}^{R_0} = 1.$ Also, as $\rho$ increases in the interval $(0,1)$, $\Upsilon_\rho^{R_0}$ is constant, $\Upsilon_{\beta_1}^{R_0}$ increases, and $\Upsilon_{\beta_2}^{R_0}$ decreases.  Thus, for larger values of $\rho$, $R_0$ becomes more sensitive to changes in $\beta_1$ and less sensitive to changes in $\beta_2$.  We can interpret this to mean controlling $R_0$ with large initial population of $S_1$ will benefit more from controlling the behavior of individuals in $S_1$ rather than those of $S_2$.

Since the structure of each sensitivity index is the same, determining with respect to which sociological parameter $R_0$ is most sensitive is determining the order of $\beta_2$, $(1-\rho)\beta_2$, and $\rho\beta_1$.  Since $1-\rho < 1$, it is always the case that $(1-\rho)\beta_2 < \beta_2$.  Thus, $\Upsilon_{\rho}^{R_0} < \Upsilon_{\beta_1}^{R_0}$, which is interpreted as $R_0$ is always more sensitive to changes in $\beta_1$ than those of $\rho$.  Further orderings of $\Upsilon_{\rho}^{R_0}, \Upsilon_{\beta_1}^{R_0},$ and $\Upsilon_{\beta_2}^{R_0}$ are dependent on parameter values.  While the next result demonstrates how the sensitivity indices are ordered, the exact values and relative magnitudes depend on exact values of $\rho, \beta_1$, and $\beta_2$. 

\begin{theorem}\label{Theorem:SensitivityBaseModel} For System \eqref{System:BasicModel}, the three sensitivity indices are ordered in the following way:
	\begin{enumerate}[font=\normalfont]
		\item[\textrm{(a)}] $\Upsilon_{\rho}^{R_0} < \Upsilon_{\beta_1}^{R_0} < \Upsilon_{\beta_2}^{R_0}$ if $\rho < \frac{\beta_2}{\beta_1+\beta_2}$.
		\item[(b)] $\Upsilon_{\rho}^{R_0} < \Upsilon_{\beta_2}^{R_0} < \Upsilon_{\beta_1}^{R_0}$ if  $\frac{\beta_2}{\beta_1+\beta_2} < \rho < \frac{\beta_2}{\beta_1}$.
		\item[(c)] $\Upsilon_{\beta_2}^{R_0} < \Upsilon_{\rho}^{R_0} < \Upsilon_{\beta_1}^{R_0}$ if $\rho > \frac{\beta_2}{\beta_1}$.
	\end{enumerate}
\end{theorem}
Once again, we can see System \eqref{System:BasicModel} undergoes a bifurcation as the parameter $\rho$ varies from 0 to 1.  In this case, the change in the sensitivity of $R_0$ to sociological parameters changes as $\rho$ increases. 

\section{A More Robust Sociological-SIR Model}\label{Section:ComplexModel}

One considerable disadvantage to the model presented in the previous section is the behavior of asymptomatic individuals.  The previous model was devised in such a way that healthy individuals were divided into two separate populations with distinct sociological behaviors that influenced the contact rates and thus the infection rates of the disease.  However, an asymptomatic person from either group was designed to behave the same regardless of the healthy population from which they originated.  While this simplified the model greatly, this does not lend to the usefulness in studying the effects of the sociological phenomenon on the spread of a disease.  To overcome this shortcoming, we devise a more robust SIR model with sociological phenomenon in this section, which we refer to by \eqref{System:ComplexModel}.

In this new model, we still have two sub-populations $S_1$ and $S_2$ of healthy individuals.  The groups $S_1$ and $S_2$ are biologically identical, in that individuals in either group are equally likely to contract the disease and recover identically.  However, the behavior of each group determines differences in infection rates.  The behavior of $S_1$ is such that the incident of infection is greater than that of $S_2$.  These rates are denoted by $\beta_1$ and $\beta_2$.  To model different types of group behavior, we assume $\beta_1 > \beta_2$, so the two groups are distinguishable.  

In this model, we have three infected sub-populations, those infected that exhibit symptoms, denoted by $I_s$ and those who are asymptomatic.  Specifically, $A_i$ is the asympomatic individuals that originated from healthy population $S_i$ (for $i = 1,2$).   In this model, any infected person is infectious.  Any infectious individual can contribute to the spread of the disease, and thus the asymptomatic individuals are denoted by $A = A_1 + A_2$ and the infectious individuals are denoted by $I = A_1+A_2+I_s$.  Since the healthy populations $S_1$ and $S_2$ are biologically identical, the proportion of infections that result in symptomatic individuals is denoted by $\lambda$, where the proportion of infections resulting in an asymptomatic individual is $1-\lambda$.  In this model, asympomatic individuals have the same rate of infection as their healthy counterparts, the assumption is sociologically individuals in $S_i$ and $A_i$ are identical.  Individuals in a healthy population can decide to change their behavior and move to the other healthy population.  We assume the rate of change between healthy populations to be denoted by $\alpha_1$ and $\alpha_2$.  As individuals in the asymptomatic classes act the same as those in the healthy classes, the same transition between $A_1$ and $A_2$ exist as for $S_1$ and $S_2$.  As with the parameters $\beta_1$ and $\beta_2$, we assume that individuals in $S_2$ (or $A_2$) will not likely change behavior, where as individuals in $S_1$ (and $A_1$) are more likely to change behavior.  Thus, $\alpha_1 > \alpha_2$.

This marks a distinct difference between this model and an SIR model where individuals in different healthy classes are biologically different.  What might be seen is a vaccine administered, which would allow for transition from one class to another, but not back (in the example of a vaccine this reverse transition would not occur on small time scales).  The interested reader is refered to \cite{MengChen:2008} as an example of such a model.

The model we propose in this section allows for recovery, which we denote by the group $R$.  We assume the rate of recovery is the same for individuals in either $A_1, A_2,$ or $I_s$, which we denote by $\kappa$.  Likewise, asymptomatic individuals develop symptoms at a rate of $\gamma$.   Finally, we assume our population to be closed with a total population $N = S_1 + S_2 + A_1 + A_2 + I_s + R$. The component diagram for \eqref{System:ComplexModel} is provided in Figure \ref{Figure:CompartmentDiagram}, and the system of ordinary differential equations is given below.

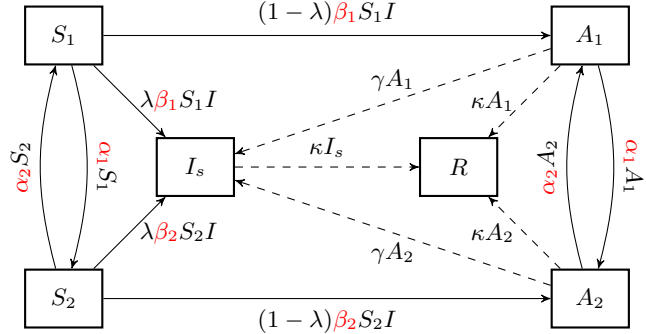
\begin{figure}
	\begin{center}
		\begin{tikzpicture} [->,>=stealth',every node/.style={font=\footnotesize}]
			\coordinate (center);
			\node[state] (Is) [left of=center] {$I_s$};
			\node[state] (R) [right of=center] {$R$};
			\node[state] (S1) [above of=Is, left of=Is] {$S_1$};
			\node[state] (S2) [below of=Is, left of=Is] {$S_2$};
			\node[state] (A1) [above of=R, right of=R] {$A_1$};
			\node[state] (A2) [below of=R, right of=R] {$A_2$};
			
			\path [bend left=15] (S1) edge node [above,rotate=-90] {$\textcolor{red}{\alpha_1} S_1$} (S2);
			\path [bend left=15] (S2) edge node [above,rotate=90] {$\textcolor{red}{\alpha_2} S_2$} (S1);
			\path [bend left=15] (A1) edge node [above,rotate=-90] {$\textcolor{red}{\alpha_1} A_1$} (A2);
			\path [bend left=15] (A2) edge node [above,rotate=90] {$\textcolor{red}{\alpha_2} A_2$} (A1);
			\path (S1) edge node [above] {$(1-\lambda)\textcolor{red}{\beta_1}S_1I$} (A1);
			\path (S1) edge node [right] {$\lambda\textcolor{red}{\beta_1}S_1I$} (Is);
			\path (S2) edge node [below] {$(1-\lambda)\textcolor{red}{\beta_2}S_2I$} (A2);
			\path (S2) edge node [right] {$\lambda\textcolor{red}{\beta_2}S_2I$} (Is);
			\path (A1) [dashed] edge node [left] {$\kappa A_1$} (R);
			\path (A2) [dashed] edge node [left] {$\kappa A_2$} (R);
			\path (A1) [dashed] edge node [above] {$\gamma A_1$} (Is);
			\path (A2) [dashed] edge node [below] {$\gamma A_2$} (Is);
			\path (Is) [dashed] edge node [above] {$\kappa I_s$} (R);
		\end{tikzpicture}
	\end{center}
	\caption{Compartment diagram of the interactions between sub-classes of the population.}\label{Figure:CompartmentDiagram}
\end{figure}

{\everymath={\displaystyle}
	\begin{equation}\label{System:ComplexModel}\tag{MB}
		\begin{array}{c@{}l}
			\frac{dS_1}{dt} &{}= \alpha_2\frac{S_2}{N} - \left[\alpha_1 + \beta_1 I\right]\frac{S_1}{N}\\[1em]
			\frac{dS_2}{dt} &{}= \alpha_1\frac{S_1}{N} - \left[\alpha_2 + \beta_2 I\right]\frac{S_2}{N}\\[1em]
			\frac{dA_1}{dt} &{}= (1-\lambda)\beta_1 I\frac{S_1}{N} + \alpha_2A_2 - (\alpha_1+\gamma+\kappa)A_1\\[1em]
			\frac{dA_2}{dt} &{}= (1-\lambda)\beta_2 I\frac{S_2}{N} + \alpha_1A_1 - (\alpha_2+\gamma+\kappa)A_2\\[1em]
			\frac{dI_s}{dt} &{}= \lambda\left[\beta_1\frac{S_1}{N} + \beta_2\frac{S_2}{N}\right]I + \gamma A - \kappa I_s\\[1em]
			\frac{dR}{dt} &{}= \kappa I
		\end{array}
\end{equation}}

Appendix B provides several numerical simulations of model \eqref{System:ComplexModel}, with an eye toward the effect of changes in the system parameters $\beta_1, \beta_2, \alpha_1,$ and $\alpha_2$.  Specifically, the reader is directed to Figures \ref{Figure:ComplexModelSimulations}, \ref{Figure:ComplexModelEffectOfB1}, \ref{Figure:ComplexModelEffectOfB2}, \ref{Figure:ComplexModelEffectOfA1}, and \ref{Figure:ComplexModelEffectOfA2}.

\subsection{Stability Analysis}
For System \eqref{System:ComplexModel} to be in equilibrium, it must be the case that $A_1 = A_2 = I_s = 0$ and $\alpha_1S_1 = \alpha_2S_2$.  Thus, System \eqref{System:ComplexModel} only exhibits disease-free equilibria, which have the form  $$(S_1^*,S_2^*,A_1^*,A_2^*,I_s^*,R^*) = (S_1^*, S_2^*, 0, 0, 0, 0),$$ with $N = S_1^* + S_2^*$ and $\alpha_1S_1 = \alpha_2S_2$.  

\begin{theorem}\label{Theorem:EquilibrimBasicModel} Equilibra of System \eqref{System:ComplexModel} have the form $(S_1=\rho N, S_2=(1-\rho)N, A_1=0, A_2=0, I_s=0, R=0)$ where \[\rho = \frac{\alpha_2}{\alpha_1+\alpha_2}.\]
\end{theorem}

The linearized infection subsystem about the disease-free equilibrium from Theorem \ref{Theorem:EquilibrimBasicModel} is given by 
{\everymath={\displaystyle}
	\[
	\begin{array}{c@{}l}
		\frac{dA_1}{dt} &{}= (1-\lambda)\beta_1 I\rho+\alpha_2A_2-(\alpha_1+\gamma+\kappa)A_1 \\[1em]
		\frac{dA_2}{dt} &{}= (1-\lambda)\beta_2 I(1-\rho)+\alpha_1A_1-(\alpha_2+\gamma+\kappa)A_2 \\[1em]
		\frac{dI_s}{dt} &{}= \lambda[\beta_1\rho+\beta_2(1-\rho)]I+\gamma A-\kappa I_s
	\end{array}
	\]}
and the Jacobian of the linearized infection subsystem by 
\[
J= \left(\begin{matrix} J_1 & J_2 & J_3\end{matrix}\right)\] where
\[\begin{aligned}
J_1 &= \left(\begin{matrix}
	(1-\lambda)\beta_1\rho-(\alpha_1+\gamma+\kappa) \\[.5em]
	(1-\lambda)\beta_2(1-\rho)+\alpha_1  \\[.5em]
	\lambda[\beta_1\rho+\beta_2(1-\rho)]+\gamma 
\end{matrix}\right),\\[.5em]]
J_2 &= \left(\begin{matrix}
	(1-\lambda)\beta_1\rho+\alpha_2  \\[.5em]
	(1-\lambda)\beta_2(1-\rho)-(\alpha_2+\gamma+\kappa) \\[.5em]
	\lambda[\beta_1\rho+\beta_2(1-\rho)]+\gamma \end{matrix}\right),\\[.5em]
J_3 &= 
\left(\begin{matrix}
	(1-\lambda)\beta_1\rho  \\[.5em]
	(1-\lambda)\beta_2(1-\rho)  \\[.5em]
	\lambda[\beta_1\rho+\beta_2(1-\rho)]-\kappa\end{matrix}\right).
\end{aligned}\] We decompose the Jacobian matrix as $J=T+\Sigma$, where $T$ is the new infection matrix and $\Sigma$ is the tranistion matrix.  Thus for $1 \leq i,j \leq 3$, $T = (T_{ij})$ where
\[T_{ij} = \begin{cases}
	(1-\lambda)\beta_1\rho & \text{if $i = 1$},\\
	(1-\lambda)\beta_2(1-\rho) & \text{if $i=2$},\\
	\lambda B_\rho & \text{if $i = 3$}
\end{cases}\]
and
\[
\Sigma =
\left(\begin{matrix}
	-(\alpha_1+\gamma+\kappa) & \alpha_2 & 0 \\[.5em]
	\alpha_1 & -(\alpha_2+\gamma+\kappa) & 0  \\[.5em]
	\gamma & \gamma &
	-\kappa\end{matrix}\right).
\]
We then form the Next Generation Matrix $K=-T\Sigma^{-1}$ and find that the eigenvalues of $K$ are 0 (of multiplicity two) and 
\[
\ell=\frac{1}{\kappa}(\beta_1\rho+\beta_2(1-\rho))
\]

\begin{theorem}\label{Theorem:ComplexModel|R0}
	For system \eqref{System:ComplexModel}, the basic reproduction number is
	\[R_0 = \frac{1}{\kappa}\left[\rho\beta_1 + (1-\rho)\beta_2\right].\] Moreover, the disease-free equilibrium \[(\rho N, (1-\rho)N, 0, 0, 0, 0)\]
	for $\rho = \displaystyle\frac{\alpha_2}{\alpha_1+\alpha_2}$ is asmptotically stable if and only if \begin{equation}\label{Equation:ComplexModel|R0ConvexityCondition}
		\rho\beta_1 + (1-\rho)\beta_2 < \kappa.
	\end{equation}
\end{theorem}

At first glance, Theorems \ref{Theorem:BaseModel:R0} and \ref{Theorem:ComplexModel|R0} look quite similar.  This is to be expected as the form of $R_0$ is the same for both systems.  One important difference to note at this point is that in System \eqref{System:BasicModel} the parameter $\rho$ could take on any value in $(0,1)$ freely.  However, in System \eqref{System:ComplexModel}, \begin{equation}\label{Inequality:UpperBoundOfRho}\rho = \frac{\alpha_2}{\alpha_1+\alpha_2} < \frac{1}{2}\end{equation} since $\alpha_1 > \alpha_2$ by assumption.

Since the form of $R_0$ is the same for the two models, much of the stability analysis holds the same for both.  The justification for Theorems \ref{Theorem:BasicModel|StabilityCriterion1} and \ref{Theorem:BasicModel|StabilityCriterion2} are the same for System \eqref{System:ComplexModel}.  So we obtain analogous results, which we combine into the following.

\begin{theorem}
	The disease-free equilibrium $(\rho N, (1-\rho)N, 0, 0, 0, 0)$ of System \eqref{System:ComplexModel} is asymptotically unstable if $\beta_2 \geq \kappa$.  Moreover, if $\beta_2 < \min\{\beta_1,\kappa\}$, then the disease-free equilibrium is asymptotically stable if and only if \[\rho < \frac{\kappa-\beta_2}{\beta_1-\beta_2},\] which is equivalent to 
	\[\alpha_1 > -\alpha_2\frac{\kappa - \beta_1}{\kappa-\beta_2}.\]
\end{theorem}

Thus, the diagrams depicted in Figure \ref{Figure:BaseModel:ActualFeasibleSetsKLarge} apply to System \eqref{System:ComplexModel} exactly the same, and so we have the following that holds.

\begin{theorem} For the System \eqref{System:ComplexModel}, if $\kappa \geq 1$ then the disease free equilibrium of the form $(\rho N, (1-\rho)N, 0, 0, 0, 0)$ is asymptotically stable for every triple $(\rho,\beta_1,\beta_2) \in \left(0,\frac12\right) \times (0,1] \times (0,1)$ with $\beta_2 < \beta_1$ and $\rho = \displaystyle\frac{\alpha_2}{\alpha_1+\alpha_2}$.
\end{theorem}

When considering the feasibile sets for System \eqref{System:ComplexModel}, the upper bound of $\rho$ established in inequality \eqref{Inequality:UpperBoundOfRho} plays a roll in shaping the possible Type of $\rho$-feasible sets.  In the case that $\kappa = \frac12$, in System \eqref{System:ComplexModel} it is not possible for $\rho \geq \kappa$.  Since it must be the case that $\rho < \kappa$, the only $\rho$-feasible set for this system is of Type 1, depicted in Figure \ref{Figure:ComplexModel|ActualFeasibleSetK1/2}. 

\begin{figure}[H]
	\begin{center}
		\begin{tikzpicture}
			\fill [gray!20] (0,0) -- (1.25, 1.25) -- (2.5,.5) -- (2.5,0) -- (0,0) -- cycle;
			\draw [dashed, thick] (0,0) -- (0,2.5);
			\draw [thick] (0,2.5) -- (0,2.75);
			\draw node [above] at (0,2.75) {\scriptsize$\beta_2$};
			\draw [dashed, thick] (0,0) -- (2.5,0);
			\draw [thick] (2.5,0) -- (2.75,0);
			\draw node [right] at (2.75,0) {\scriptsize$\beta_1$};
			\draw [thick] (2.5,.07) -- (2.5,-.07) node [below] at (2.5,0) {\tiny$1$};
			\draw [thick] (.07,2.5) -- (-.07,2.5) node [left] at (0,2.5) {\tiny$1$};
			\draw [dashed, thick] (0,2.5) -- (2.5,2.5);
			\draw [thick] (2.5,0) -- (2.5,2.5);
			\draw [dashed, thick, blue] (0,0) -- (2.5,2.5);
			\draw [dashed, thick, red] (0,2) -- (2.5,.5);
			\draw [thick] (-.07,2) -- (.07,2) node [left] at (-.1,2) {\tiny$\frac{\kappa}{1-\rho}$};
			\draw [thick] (2.43,.5) -- (2.57,.5)  node [right] at (2.5,.7) {\tiny$\frac{\kappa-\rho}{1-\rho}$};
			\draw [dotted] (0,1.25) -- (1.25,1.25);
			\draw [thick] (-.07,1.25) -- (.07,1.25) node [left] at (0,1.25) {\tiny$\kappa$};
			\draw node [below] at (2,0) {\textcolor{white}{\tiny$\frac{\kappa}{\rho}$}};
			\draw node [above] at (.85,2.5) {\textcolor{white}{\tiny$\frac{\kappa-(1-\rho)}{\rho}$}};
		\end{tikzpicture}	
		\caption{Feasible set of System \eqref{System:ComplexModel} for $\kappa = \frac12$.}\label{Figure:ComplexModel|ActualFeasibleSetK1/2}
	\end{center}
\end{figure}

The feasible sets for the cases of $0 < \kappa < \frac12$ and $\frac12 < \kappa < 1$ are shown in Figures \ref{Figure:ComplexModel:ActualFeasibleSetsK<1/2} and \ref{Figure:ComplexModel:ActualFeasibleSetsK>1/2}, respectively.  Again we see the upper bound of $\rho$ plays a role in eliminating possible feasible sets that were exhibited in System \eqref{System:BasicModel}.  

The bifurcations in the feasible sets of System \eqref{System:ComplexModel} also changes as a result in the bound on $\rho$ as well as the elimination of possible feasible sets.  In the case of $\kappa = \frac12$, there is no bifurcation in feasible set type when the value of $\rho$ varies in $\left(0,\frac12\right)$, unlike in System \ref{System:BasicModel}.  The bifurcation diagrams for the cases of $0 < \kappa < \frac12$ and $\frac12 < \kappa < 1$ are shown in Figure \ref{Figure:BifurcationComplexModel}.

\subsection{Sensitivity Analysis}

For System \eqref{System:ComplexModel}, we wish to determine the sensitivity of $R_0$ to the sociological parameters $\rho, \beta_1,$ and $\beta_2$, as we did for System \eqref{System:BasicModel}.  Since the form of $R_0$ is the same for both systems, the sensitivity indices $\Upsilon_\rho^{R_0}, \Upsilon_{\beta_1}^{R_0},$ and $\Upsilon_{\beta_2}^{R_0}$ are the same as those calculated in \eqref{Equation:MASensitivityIndices}.  

However, unlike in the previous system, System \eqref{System:ComplexModel} has two more sociological parameters, namely $\alpha_1$ and $\alpha_2$, of which $\rho$ is a function.  Thus, we wish to also compute the sensitivity indices of $R_0$ with respect to $\alpha_1$ and $\alpha_2$, and determine a relative order of sensitivity.  We calculate $\Upsilon_{\alpha_1}^{R_0}$ and $\Upsilon_{\alpha_2}^{R_0}$ thusly

\[\begin{aligned}
\Upsilon_{\alpha_1}^{R_0} &= \left(\frac{\alpha_1}{R_0}\right)\frac{\partial R_0}{\partial\alpha_1} = \frac{-\rho(1-\rho)(\beta_1-\beta_2)}{R_0\kappa}\\
\Upsilon_{\alpha_2}^{R_0} &= \left(\frac{\alpha_2}{R_0}\right)\frac{\partial R_0}{\partial\alpha_2} = \frac{\rho(1-\rho)(\beta_1-\beta_2)}{R_0\kappa}.
\end{aligned}\]

It is immediate that $\Upsilon_{\alpha_1}^{R_0} < \Upsilon_{\alpha_2}^{R_0} < \Upsilon_\rho^{R_0}.$ Thus, we see that $R_0$ will always be the least sensitive to $\alpha_1$ since $\Upsilon_{\alpha_1}^{R_0}$ is negative while the other indicies are strictly positive.  The order relation of $\Upsilon_\rho^{R_0}, \Upsilon_{\beta_1}^{R_0},$ and $\Upsilon_{\beta_2}^{R_0}$ from Theorem \ref{Theorem:SensitivityBaseModel} still hold.  So we now determine how $\Upsilon_{\alpha_2}^{R_0}$ is ordered.

\begin{theorem}\label{Theorem:SensitivityComplexSystem} For System \eqref{System:ComplexModel}, the five sensitivity indices are ordered in the following way:
	\begin{enumerate}[font=\normalfont]
		\item[(a)] $\Upsilon_{\alpha_1}^{R_0} < \Upsilon_{\alpha_2}^{R_0} < \Upsilon_{\rho}^{R_0} < \Upsilon_{\beta_1}^{R_0} < \Upsilon_{\beta_2}^{R_0}$ if $\rho < \frac{\beta_2}{\beta_1+\beta_2}$.
		\item[(b)] $\Upsilon_{\alpha_1}^{R_0} < \Upsilon_{\alpha_2}^{R_0} < \Upsilon_{\rho}^{R_0} < \Upsilon_{\beta_2}^{R_0} < \Upsilon_{\beta_1}^{R_0}$ if $\frac{\beta_2}{\beta_1+\beta_2} < \rho < \frac{\beta_2}{\beta_1}$.
		\item[(c)] $\Upsilon_{\alpha_1}^{R_0} < \Upsilon_{\alpha_2}^{R_0} < \Upsilon_{\beta_2}^{R_0} < \Upsilon_{\rho}^{R_0} <   \Upsilon_{\beta_1}^{R_0}$ if $\frac{\beta_2}{\beta_1} < \rho < \frac{\beta_2}{\beta_1-\beta_2}$.
		\item[(d)] $\Upsilon_{\alpha_1}^{R_0} < \Upsilon_{\beta_2}^{R_0} < \Upsilon_{\alpha_2}^{R_0} < \Upsilon_{\rho}^{R_0} <  \Upsilon_{\beta_1}^{R_0}$ if $\rho > \frac{\beta_2}{\beta_1-\beta_2}$.
	\end{enumerate}
\end{theorem}

Unlike in Model \eqref{System:BasicModel}, the value of $\rho$ for System \eqref{System:ComplexModel} is bounded above by $\frac12$.  Thus certain orderings of the sensitivity indices will not be possible, depending on the values of $\beta_1$ and $\beta_2$.  Specifically, if $\beta_2 \geq \frac{1}{3}\beta_1$, then case (d) of Theorem \ref{Theorem:SensitivityComplexSystem} is not possible.  Furthermore, if $\beta_2 \geq \frac{1}{2}\beta_1$, then case (c) of Theorem \ref{Theorem:SensitivityComplexSystem} is not possible as well.

\section{Conclusion}\label{Section:Conclusion}
The models proposed in this paper can be used to determine effects of sociological interventions on the propagation of a disease.  These include the implementation of ``stay in place" orders, mask mandates, and even vaccine mandates.  In \cite{Payne:2020}, the COVID-19 infection aboard the USS Theodore Roosevelt in April 2020 was studied, specifically determining the efficacy of social interventions (mask wearing, avoiding common areas, and observing social distancing) on the infection rate, as outlined in Table \ref{Table:Covid19-Rates}.

\begin{table}[H]
\caption{Infection rates associated to three sociological mitigation strategies, make wearing, avoiding common areas, and social distancing, observed on the USS Theodore Roosevelt during a COVID-19 outbreak in April 2020 \cite{Payne:2020}.}\label{Table:Covid19-Rates}
\begin{center}
\begin{tabular}{r c c}
& \multicolumn{2}{c}{\small{Infection Rates}}\\\cline{2-3}
& \footnotesize{not implemented} & \footnotesize{implemented}\\\cline{2-3}
\footnotesize{mask wearing} & \footnotesize{80.8\%} & \footnotesize{55.8\%}\\
\footnotesize{avoiding common areas} & \footnotesize{67.5\%} & \footnotesize{53.8\%}\\
\footnotesize{social distancing} & \footnotesize{70.0\%} & \footnotesize{54.7\%}\\
\end{tabular}
\end{center}
\end{table}

\begin{figure}[H]
\begin{picbox}
    \begin{center}
	\includegraphics[width=\textwidth]{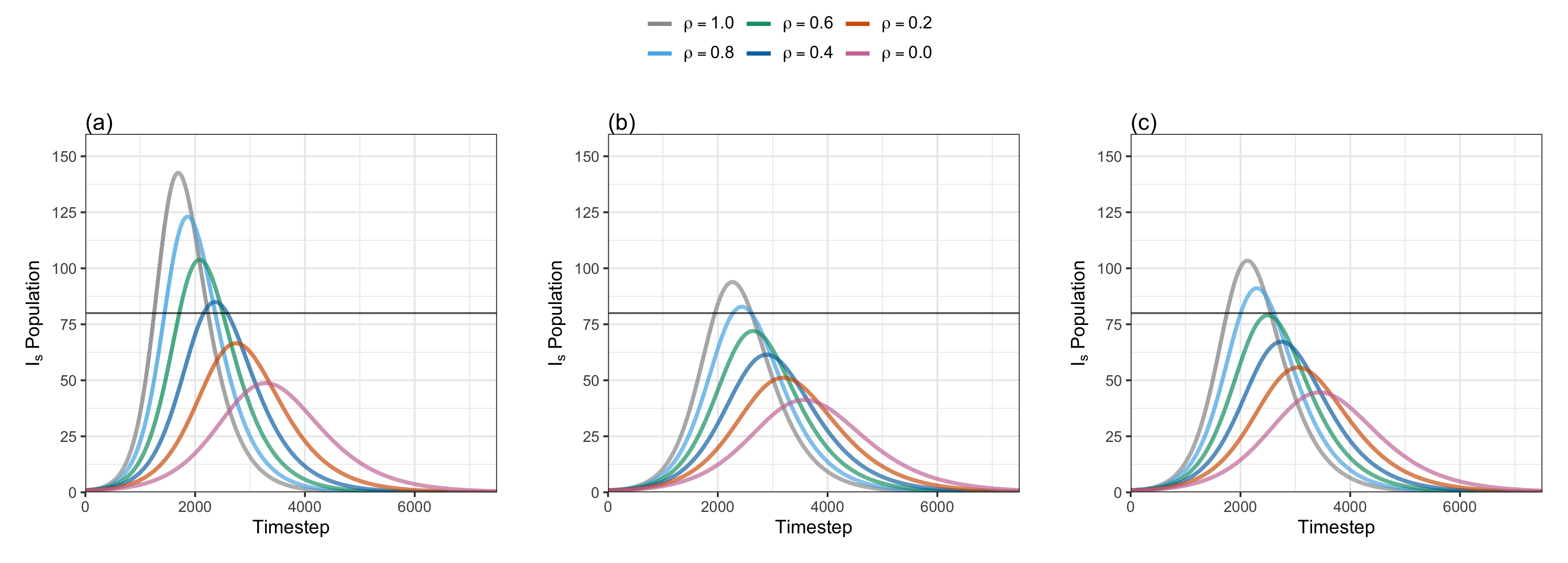}
	\end{center}
	\tcblower
	\noindent\begin{tabular}{lx{.12\textwidth}x{.12\textwidth}x{.1\textwidth}x{.1\textwidth}x{.1\textwidth}x{.1\textwidth}x{.1\textwidth}}
	    (a) & $\beta_1$={\footnotesize 0.00808} & $\beta_2$={\footnotesize 0.00558} & $\alpha_1$={\footnotesize 0.0001} & $\alpha_2$={\footnotesize 0.0} & $\gamma$={\footnotesize 0.0001} & $\lambda$={\footnotesize 0.65} & $\kappa$={\footnotesize 0.0002}
	    \\
        (b) & $\beta_1$={\footnotesize 0.00675} & $\beta_2$={\footnotesize 0.00538} & $\alpha_1$={\footnotesize 0.0001} & $\alpha_2$={\footnotesize 0.0} & $\gamma$={\footnotesize 0.0001} & $\lambda$={\footnotesize 0.65} & $\kappa$={\footnotesize 0.0002}
        \\
        (c) & $\beta_1$={\footnotesize 0.00700} & $\beta_2$={\footnotesize 0.00547} & $\alpha_1$={\footnotesize 0.0001} & $\alpha_2$={\footnotesize 0.0} & $\gamma$={\footnotesize 0.0001} & $\lambda$={\footnotesize 0.65} & $\kappa$={\footnotesize 0.0002}
	\end{tabular}
\end{picbox}
\caption{Simulations of model \eqref{System:ComplexModel} subject to the infection rates of (a) mask wearing, (b) avoiding common areas, and (c) social distancing from Table \ref{Table:Covid19-Rates}.}\label{Figure:COVID-19Parameters}
\end{figure}

This data is a means by which to estimate the parameters $\beta_1$ and $\beta_2$ in both \eqref{System:BasicModel} and \eqref{System:ComplexModel}.  In Figure \ref{Figure:COVID-19Parameters}, we perform simulations of \eqref{System:ComplexModel} to model the spread of a hypothetical disease under the three sociological mitigation strategies.  In these simulations, $S_1$ represents individuals who do not implement the mitigation strategy, versus the population $S_2$ of individuals that do implement the strategy.  While it is unclear in the original study whether these infection rates are independent of each other, we use the values as a means of illustrating how the models in this manuscript can be used to determine the most appropriate mitigation strategy.  

We devise the following scenario as means to show how the models \eqref{System:BasicModel} and \eqref{System:ComplexModel} can be utilized by policy makers to determine which sociological mitigation is most appropriate.  Suppose a closed population is subject to the disease modelled by \eqref{System:ComplexModel}, and that any symptomatic individual will require medical resources.  Also suppose for this discussion that the medical resources can support at most 80 symptomatic individuals at any one time.  We can see from Figure \ref{Figure:COVID-19Parameters}, each mitigation strategy can result in the peak of the symptomatic population being under the threshold value of 80.  The difference between the mitigation strategies is the level of participation required to ``flatten the curve."  From Figure \ref{Figure:COVID-19Parameters}, more than 60\% of the population will be required to wear masks for the infection curve to peak under the threshold.  In comparison, it will require only about 20\% of the population to either avoiding common areas or practice social distancing to have the same effect.

So, to decide which strategy to implement, understanding of the cost of implementation is required.  In the example of the USS Theodore Roosevelt, it may be easier to implement a mask requirement rather than mandate social distancing or avoiding common areas on a ship, where physically this may be impossible.  So even through it would require less of the population to implement social distancing to flatten the infection curve enough to not stress medical resources, it may be infeasible or impossible.  Thus, the appropriate mitigation strategy, and the level of participation, must be determined by those who can determine the cost-benefit of each strategy.  

\section{Future Work}\label{Section:Future}
The use of the models developed in this paper only help explore situations once the disease is present in the population.  A more realistic approach might be to first model the population with a single susceptible population, the results of which can then be used to develop the initial conditions for the sociological SIR models introduced in this paper. This would rely on numerical simulation over the basic reproduction number, as initial conditions would no longer be near equilibrium.

\begin{figure}[H]
\begin{picbox}
	\includegraphics[width=\textwidth]{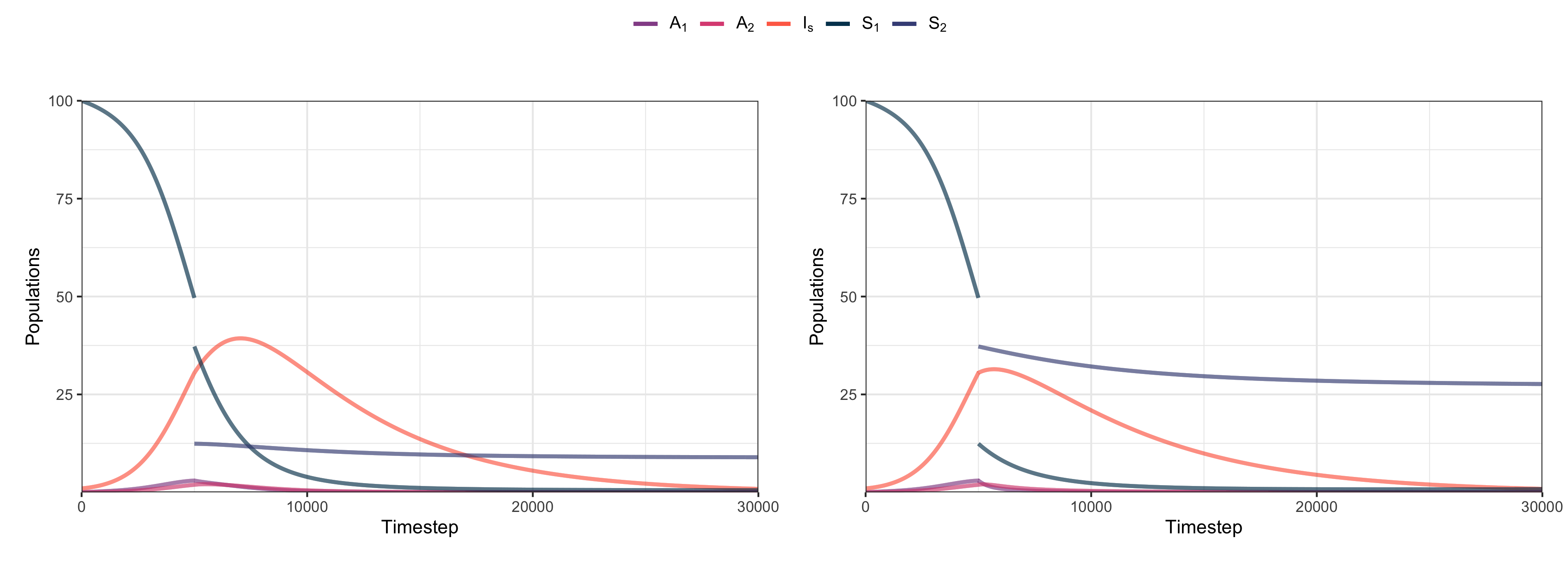}
	\tcblower
	\noindent\begin{tabular}{x{.115\textwidth}x{.115\textwidth}x{.115\textwidth}x{.115\textwidth}x{.115\textwidth}x{.115\textwidth}x{.115\textwidth}}
	    $\beta_1$={\footnotesize 0.0011} & $\beta_2$={\footnotesize 0.0001} & $\alpha_1$={\footnotesize 0.001} & $\alpha_2$={\footnotesize 0.0001} & $\gamma$={\footnotesize 0.0001} & $\lambda$={\footnotesize 0.65} & $\kappa$={\footnotesize 0.0002}
	\end{tabular}
\end{picbox}
\caption{Simulations of mixed SIR models with change in model occurring at the same time with $\rho = 0.25$ and $\rho = 0.75$.}\label{Figure:MixedModelRho}
\end{figure}

\begin{figure}[H]
\begin{picbox}
	\includegraphics[width=\textwidth]{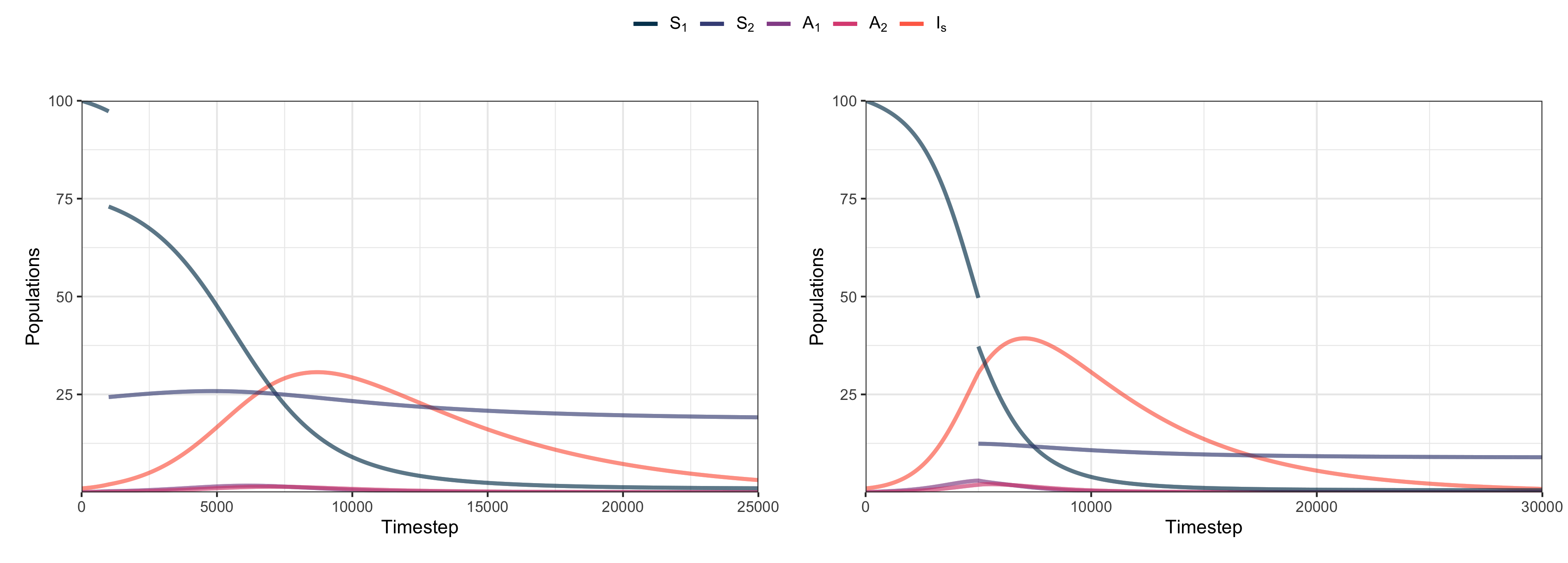}
	\tcblower
	\noindent\begin{tabular}{x{.115\textwidth}x{.115\textwidth}x{.115\textwidth}x{.115\textwidth}x{.115\textwidth}x{.115\textwidth}x{.115\textwidth}}
	    $\beta_1$={\footnotesize 0.0011} & $\beta_2$={\footnotesize 0.0001} & $\alpha_1$={\footnotesize 0.001} & $\alpha_2$={\footnotesize 0.0001} & $\gamma$={\footnotesize 0.0001} & $\lambda$={\footnotesize 0.65} & $\kappa$={\footnotesize 0.0002}
	\end{tabular}
\end{picbox}
\caption{Simulations of mixed SIR models with change in model occurring at timesteps 1000 and 5000 ($\rho = 0.25$).}\label{Figure:MixedModelTime}
\end{figure}

As an example, we first simulate the disease propagation through the population with $\rho=1$, that is a single susceptible-population model. We then take the resulting population profile as the initial condition for the Sociological SIR model described in Section \ref{Section:ComplexModel}.  
In Figure \ref{Figure:MixedModelRho} we simulate changing the value of $\rho$ from $0.25$ to $0.75$.  In Figure \ref{Figure:MixedModelTime} we simulate a change when the system switches from one susceptible population to two. 

This type of mixed model can be used to determine the sociological response necessary depending on when the disease is established in the population, as well as the magnitude of the control mechanisms to ``flatten the curve", for example.

\bibliographystyle{amsplain}
\bibliography{references.bib}

\newpage
\section*{Appendix A: Feasibility Sets and Bifurcation Diagrams}

\newpage
\begin{figure}[H]
	\begin{center}
		\begin{tabular}{rc|cl}
			\multicolumn{1}{c}{} & \multicolumn{1}{c}{\small$0 < \kappa < \frac12$} & \multicolumn{1}{c}{\small$\frac12 < \kappa < 1$} & \multicolumn{1}{c}{}\\[3mm]
			\hline
			{\small$0 < \rho < \kappa$}
			&
			\adjustbox{valign=c}{%
				\begin{tikzpicture}
					\fill [gray!20] (0,0) -- (1.25, 1.25) -- (2.5,.5) -- (2.5,0) -- (0,0) -- cycle;
					\draw [dashed, thick] (0,0) -- (0,2.5);
					\draw [thick] (0,2.5) -- (0,2.75);
					\draw node [above] at (0,2.75) {\scriptsize$\beta_2$};
					\draw [dashed, thick] (0,0) -- (2.5,0);
					\draw [thick] (2.5,0) -- (2.75,0);
					\draw node [right] at (2.75,0) {\scriptsize$\beta_1$};
					\draw [thick] (2.5,.07) -- (2.5,-.07) node [below] at (2.5,0) {\scriptsize$1$};
					\draw [thick] (.07,2.5) -- (-.07,2.5) node [left] at (0,2.5) {\scriptsize$1$};
					\draw [dashed, thick] (0,2.5) -- (2.5,2.5);
					\draw [thick] (2.5,0) -- (2.5,2.5);
					\draw [dashed, thick, blue] (0,0) -- (2.5,2.5);
					\draw [dashed, thick, red] (0,2) -- (2.5,.5);
					\draw [thick] (-.07,2) -- (.07,2) node [left] at (-.1,2) {\scriptsize$\frac{\kappa}{1-\rho}$};
					\draw [thick] (2.43,.5) -- (2.57,.5)  node [right] at (2.5,.7) {\scriptsize$\frac{\kappa-\rho}{1-\rho}$};
					\draw [dotted] (0,1.25) -- (1.25,1.25);
					\draw [thick] (-.07,1.25) -- (.07,1.25) node [left] at (0,1.25) {\scriptsize$\kappa$};
					\draw node [below] at (2,0) {\textcolor{white}{\scriptsize$\kappa/\rho$}};
					\draw node [above] at (1,2.5) {\textcolor{white}{\scriptsize$\frac{\kappa-(1-\rho)}{\rho}$}};	
				\end{tikzpicture}
			}
			&
			\adjustbox{valign=c}{%
				\begin{tikzpicture}
					\fill [gray!20] (0,0) -- (1.25, 1.25) -- (2,0) -- (0,0) -- cycle;
					\draw [dashed, thick] (0,0) -- (0,2.5);
					\draw [thick] (0,2.5) -- (0,2.75);
					\draw node [above] at (0,2.75) {\scriptsize$\beta_2$};
					\draw [dashed, thick] (0,0) -- (2.5,0);
					\draw [thick] (2.5,0) -- (2.75,0);
					\draw node [right] at (2.75,0) {\scriptsize$\beta_1$};
					\draw [thick] (2.5,.07) -- (2.5,-.07) node [below] at (2.5,0) {\scriptsize$1$};
					\draw [thick] (.07,2.5) -- (-.07,2.5) node [left] at (0,2.5) {\scriptsize$1$};
					\draw [dashed, thick] (0,2.5) -- (2.5,2.5);
					\draw [thick] (2.5,0) -- (2.5,2.5);
					\draw [dashed, thick, blue] (0,0) -- (2.5,2.5);
					\draw [dashed, thick, red] (.5,2.5) -- (2,0);
					\draw [thick] (.5,2.57) -- (.5,2.43) node [above] at (1,2.5) {\scriptsize$\frac{\kappa-(1-\rho)}{\rho}$};
					\draw node [left] at (-.1,2) {\textcolor{white}{\scriptsize$\frac{\kappa}{1-\rho}$}};
					\draw [thick] (2,-.07) -- (2,.07) node [below] at (2,0) {\scriptsize$\kappa/\rho$};
					\draw [dotted] (0,1.25) -- (1.25,1.25);
					\draw [thick] (-.07,1.25) -- (.07,1.25) node [left] at (0,1.25) {\scriptsize$\kappa$};
					\draw node [left] at (0,2) {\textcolor{white}{\scriptsize$\kappa/\rho$}};
					\draw node [right] at (2.5,.7) {\textcolor{white}{\scriptsize$\frac{\kappa-\rho}{1-\rho}$}};
				\end{tikzpicture}
			}
			&
			{\small$0 < \rho < 1-\kappa$}
			\\\hline	
			{\small$\rho=\kappa$}
			&
			\adjustbox{valign=c}{%
				\begin{tikzpicture}
					\fill [gray!20] (0,0) -- (1.1, 1.1) -- (2.5,0) -- (0,0) -- cycle;
					\draw [dashed, thick] (0,0) -- (0,2.5);
					\draw [thick] (0,2.5) -- (0,2.75);
					\draw node [above] at (0,2.75) {\scriptsize$\beta_2$};
					\draw [dashed, thick] (0,0) -- (2.5,0);
					\draw [thick] (2.5,0) -- (2.75,0);
					\draw node [right] at (2.75,0) {\scriptsize$\beta_1$};
					\draw [thick] (2.5,.07) -- (2.5,-.07) node [below] at (2.5,0) {\scriptsize$1$};
					\draw [thick] (.07,2.5) -- (-.07,2.5) node [left] at (0,2.5) {\scriptsize$1$};
					\draw [dashed, thick] (0,2.5) -- (2.5,2.5);
					\draw [thick] (2.5,0) -- (2.5,2.5);
					\draw [dashed, thick, blue] (0,0) -- (2.5,2.5);
					\draw [dashed, thick, red] (0,2) -- (2.5,0);
					\draw [thick] (-.07,2) -- (.07,2) node [left] at (-.1,2) {\scriptsize$\frac{\kappa}{1-\rho}$};
					\draw [dotted] (0,1.1) -- (1.1,1.1);
					\draw [thick] (-.07,1.1) -- (.07,1.1) node [left] at (0,1.1) {\scriptsize$\kappa$};
					\draw node [below] at (2,0) {\textcolor{white}{\scriptsize$\kappa/\rho$}};
					\draw node [right] at (2.5,.7) {\textcolor{white}{\scriptsize$\frac{\kappa-\rho}{1-\rho}$}};
					\draw node [above] at (1,2.5) {\textcolor{white}{\scriptsize$\frac{\kappa-(1-\rho)}{\rho}$}};	
				\end{tikzpicture}
			}
			&
			\adjustbox{valign=c}{%
				\begin{tikzpicture}
					\fill [gray!20] (0,0) -- (1.1, 1.1) -- (2,0) -- (0,0) -- cycle;
					\draw [dashed, thick] (0,0) -- (0,2.5);
					\draw [thick] (0,2.5) -- (0,2.75);
					\draw node [above] at (0,2.75) {\scriptsize$\beta_2$};
					\draw [dashed, thick] (0,0) -- (2.5,0);
					\draw [thick] (2.5,0) -- (2.75,0);
					\draw node [right] at (2.75,0) {\scriptsize$\beta_1$};
					\draw [thick] (2.5,.07) -- (2.5,-.07) node [below] at (2.5,0) {\scriptsize$1$};
					\draw [thick] (.07,2.5) -- (-.07,2.5) node [left] at (0,2.5) {\scriptsize$1$};
					\draw [dashed, thick] (0,2.5) -- (2.5,2.5);
					\draw [thick] (2.5,0) -- (2.5,2.5);
					\draw [dashed, thick, blue] (0,0) -- (2.5,2.5);
					\draw [dashed, thick, red] (0,2.5) -- (2,0);
					\draw [thick] (2,-.07) -- (2,.07) node [below] at (2,0) {\scriptsize$\kappa/\rho$};
					\draw [dotted] (0,1.1) -- (1.1,1.1);
					\draw [thick] (-.07,1.25) -- (.07,1.25) node [left] at (0,1.25) {\scriptsize$\kappa$};
					\draw node [left] at (-.1,2) {\textcolor{white}{\scriptsize$\frac{\kappa}{1-\rho}$}};
					\draw node [above] at (1,2.5) {\textcolor{white}{\scriptsize$\frac{\kappa-(1-\rho)}{\rho}$}};
					\draw node [right] at (2.5,.7) {\textcolor{white}{\scriptsize$\frac{\kappa-\rho}{1-\rho}$}};
				\end{tikzpicture}
			}
			&
			{\small$\rho = 1-\kappa$}
			\\\hline
			{\small$\kappa < \rho < 1-\kappa$}
			&
			\adjustbox{valign=c}{%
				\begin{tikzpicture}
					\fill [gray!20] (0,0) -- (1, 1) -- (2,0) -- (0,0) -- cycle;
					\draw [dashed, thick] (0,0) -- (0,2.5);
					\draw [thick] (0,2.5) -- (0,2.75);
					\draw node [above] at (0,2.75) {\scriptsize$\beta_2$};
					\draw [dashed, thick] (0,0) -- (2.5,0);
					\draw [thick] (2.5,0) -- (2.75,0);
					\draw node [right] at (2.75,0) {\scriptsize$\beta_1$};
					\draw [thick] (2.5,.07) -- (2.5,-.07) node [below] at (2.5,0) {\scriptsize$1$};
					\draw [thick] (.07,2.5) -- (-.07,2.5) node [left] at (0,2.5) {\scriptsize$1$};
					\draw [dashed, thick] (0,2.5) -- (2.5,2.5);
					\draw [thick] (2.5,0) -- (2.5,2.5);
					\draw [dashed, thick, blue] (0,0) -- (2.5,2.5);
					\draw [dashed, thick, red] (0,2) -- (2,0);
					\draw [thick] (-.07,2) -- (.07,2) node [left] at (-.1,2) {\scriptsize$\frac{\kappa}{1-\rho}$};
					\draw [dotted] (0,1) -- (1,1);
					\draw [thick] (-.07,1) -- (.07,1) node [left] at (0,1) {\scriptsize$\kappa$};
					\draw [thick] (2,-.07) -- (2,.07) node [below] at (2,0) {\scriptsize$\kappa/\rho$};
					\draw node [right] at (2.5,.7) {\textcolor{white}{\scriptsize$\frac{\kappa-\rho}{1-\rho}$}};
					\draw node [above] at (1,2.5) {\textcolor{white}{\scriptsize$\frac{\kappa-(1-\rho)}{\rho}$}};	
				\end{tikzpicture}
			}
			&
			\adjustbox{valign=c}{%
				\begin{tikzpicture}
					\fill [gray!20] (0,0) -- (1, 1) -- (2,0) -- (0,0) -- cycle;
					\draw [dashed, thick] (0,0) -- (0,2.5);
					\draw [thick] (0,2.5) -- (0,2.75);
					\draw node [above] at (0,2.75) {\scriptsize$\beta_2$};
					\draw [dashed, thick] (0,0) -- (2.5,0);
					\draw [thick] (2.5,0) -- (2.75,0);
					\draw node [right] at (2.75,0) {\scriptsize$\beta_1$};
					\draw [thick] (2.5,.07) -- (2.5,-.07) node [below] at (2.5,0) {\scriptsize$1$};
					\draw [thick] (.07,2.5) -- (-.07,2.5) node [left] at (0,2.5) {\scriptsize$1$};
					\draw [dashed, thick] (0,2.5) -- (2.5,2.5);
					\draw [thick] (2.5,0) -- (2.5,2.5);
					\draw [dashed, thick, blue] (0,0) -- (2.5,2.5);
					\draw [dashed, thick, red] (0,2) -- (2,0);
					\draw [thick] (-.07,2) -- (.07,2) node [left] at (-.1,2) {\scriptsize$\frac{\kappa}{1-\rho}$};
					\draw [dotted] (0,1) -- (1,1);
					\draw [thick] (-.07,1) -- (.07,1) node [left] at (0,1) {\scriptsize$\kappa$};
					\draw [thick] (2,-.07) -- (2,.07) node [below] at (2,0) {\scriptsize$\kappa/\rho$};
					\draw node [right] at (2.5,.7) {\textcolor{white}{\scriptsize$\frac{\kappa-\rho}{1-\rho}$}};
					\draw node [above] at (1,2.5) {\textcolor{white}{\scriptsize$\frac{\kappa-(1-\rho)}{\rho}$}};	
				\end{tikzpicture}
			}
			&
			{\small$1-\kappa < \rho < \kappa$}
			\\\hline
			{\small$\rho = 1-\kappa$}
			&
			\adjustbox{valign=c}{%
				\begin{tikzpicture}
					\fill [gray!20] (0,0) -- (1.1, 1.1) -- (2,0) -- (0,0) -- cycle;
					\draw [dashed, thick] (0,0) -- (0,2.5);
					\draw [thick] (0,2.5) -- (0,2.75);
					\draw node [above] at (0,2.75) {\scriptsize$\beta_2$};
					\draw [dashed, thick] (0,0) -- (2.5,0);
					\draw [thick] (2.5,0) -- (2.75,0);
					\draw node [right] at (2.75,0) {\scriptsize$\beta_1$};
					\draw [thick] (2.5,.07) -- (2.5,-.07) node [below] at (2.5,0) {\scriptsize$1$};
					\draw [thick] (.07,2.5) -- (-.07,2.5) node [left] at (0,2.5) {\scriptsize$1$};
					\draw [dashed, thick] (0,2.5) -- (2.5,2.5);
					\draw [thick] (2.5,0) -- (2.5,2.5);
					\draw [dashed, thick, blue] (0,0) -- (2.5,2.5);
					\draw [dashed, thick, red] (0,2.5) -- (2,0);
					\draw [thick] (2,-.07) -- (2,.07) node [below] at (2,0) {\scriptsize$\kappa/\rho$};
					\draw [dotted] (0,1.1) -- (1.1,1.1);
					\draw [thick] (-.07,1.25) -- (.07,1.25) node [left] at (0,1.25) {\scriptsize$\kappa$};
					\draw node [left] at (-.1,2) {\textcolor{white}{\scriptsize$\frac{\kappa}{1-\rho}$}};
					\draw node [above] at (1,2.5) {\textcolor{white}{\scriptsize$\frac{\kappa-(1-\rho)}{\rho}$}};
					\draw node [right] at (2.5,.7) {\textcolor{white}{\scriptsize$\frac{\kappa-\rho}{1-\rho}$}};
				\end{tikzpicture}
			}
			&
			\adjustbox{valign=c}{%
				\begin{tikzpicture}
					\fill [gray!20] (0,0) -- (1.1, 1.1) -- (2.5,0) -- (0,0) -- cycle;
					\draw [dashed, thick] (0,0) -- (0,2.5);
					\draw [thick] (0,2.5) -- (0,2.75);
					\draw node [above] at (0,2.75) {\scriptsize$\beta_2$};
					\draw [dashed, thick] (0,0) -- (2.5,0);
					\draw [thick] (2.5,0) -- (2.75,0);
					\draw node [right] at (2.75,0) {\scriptsize$\beta_1$};
					\draw [thick] (2.5,.07) -- (2.5,-.07) node [below] at (2.5,0) {\scriptsize$1$};
					\draw [thick] (.07,2.5) -- (-.07,2.5) node [left] at (0,2.5) {\scriptsize$1$};
					\draw [dashed, thick] (0,2.5) -- (2.5,2.5);
					\draw [thick] (2.5,0) -- (2.5,2.5);
					\draw [dashed, thick, blue] (0,0) -- (2.5,2.5);
					\draw [dashed, thick, red] (0,2) -- (2.5,0);
					\draw [thick] (-.07,2) -- (.07,2) node [left] at (-.1,2) {\scriptsize$\frac{\kappa}{1-\rho}$};
					\draw [dotted] (0,1.1) -- (1.1,1.1);
					\draw [thick] (-.07,1.1) -- (.07,1.1) node [left] at (0,1.1) {\scriptsize$\kappa$};
					\draw node [below] at (2,0) {\textcolor{white}{\scriptsize$\kappa/\rho$}};
					\draw node [right] at (2.5,.7) {\textcolor{white}{\scriptsize$\frac{\kappa-\rho}{1-\rho}$}};
					\draw node [above] at (1,2.5) {\textcolor{white}{\scriptsize$\frac{\kappa-(1-\rho)}{\rho}$}};	
				\end{tikzpicture}
			}
			&
			{\small$\rho = \kappa$}
			\\\hline
			{\small$1-\kappa < \rho < 1$}
			&
			\adjustbox{valign=c}{%
				\begin{tikzpicture}
					\fill [gray!20] (0,0) -- (1.25, 1.25) -- (2,0) -- (0,0) -- cycle;
					\draw [dashed, thick] (0,0) -- (0,2.5);
					\draw [thick] (0,2.5) -- (0,2.75);
					\draw node [above] at (0,2.75) {\scriptsize$\beta_2$};
					\draw [dashed, thick] (0,0) -- (2.5,0);
					\draw [thick] (2.5,0) -- (2.75,0);
					\draw node [right] at (2.75,0) {\scriptsize$\beta_1$};
					\draw [thick] (2.5,.07) -- (2.5,-.07) node [below] at (2.5,0) {\scriptsize$1$};
					\draw [thick] (.07,2.5) -- (-.07,2.5) node [left] at (0,2.5) {\scriptsize$1$};
					\draw [dashed, thick] (0,2.5) -- (2.5,2.5);
					\draw [thick] (2.5,0) -- (2.5,2.5);
					\draw [dashed, thick, blue] (0,0) -- (2.5,2.5);
					\draw [dashed, thick, red] (.5,2.5) -- (2,0);
					\draw [thick] (.5,2.57) -- (.5,2.43) node [above] at (1,2.5) {\scriptsize$\frac{\kappa-(1-\rho)}{\rho}$};
					\draw [thick] (2,-.07) -- (2,.07) node [below] at (2,0) {\scriptsize$\kappa/\rho$};
					\draw [dotted] (0,1.25) -- (1.25,1.25);
					\draw [thick] (-.07,1.25) -- (.07,1.25) node [left] at (0,1.25) {\scriptsize$\kappa$};
					\draw node [left] at (-.1,2) {\textcolor{white}{\scriptsize$\frac{\kappa}{1-\rho}$}};
					\draw node [right] at (2.5,.7) {\textcolor{white}{\scriptsize$\frac{\kappa-\rho}{1-\rho}$}};
				\end{tikzpicture}
			}
			&
			\adjustbox{valign=c}{%
				\begin{tikzpicture}
					\fill [gray!20] (0,0) -- (1.25, 1.25) -- (2.5,.5) -- (2.5,0) -- (0,0) -- cycle;
					\draw [dashed, thick] (0,0) -- (0,2.5);
					\draw [thick] (0,2.5) -- (0,2.75);
					\draw node [above] at (0,2.75) {\scriptsize$\beta_2$};
					\draw [dashed, thick] (0,0) -- (2.5,0);
					\draw [thick] (2.5,0) -- (2.75,0);
					\draw node [right] at (2.75,0) {\scriptsize$\beta_1$};
					\draw [thick] (2.5,.07) -- (2.5,-.07) node [below] at (2.5,0) {\scriptsize$1$};
					\draw [thick] (.07,2.5) -- (-.07,2.5) node [left] at (0,2.5) {\scriptsize$1$};
					\draw [dashed, thick] (0,2.5) -- (2.5,2.5);
					\draw [thick] (2.5,0) -- (2.5,2.5);
					\draw [dashed, thick, blue] (0,0) -- (2.5,2.5);
					\draw [dashed, thick, red] (0,2) -- (2.5,.5);
					\draw [thick] (-.07,2) -- (.07,2) node [left] at (-.1,2) {\scriptsize$\frac{\kappa}{1-\rho}$};
					\draw [thick] (2.43,.5) -- (2.57,.5)  node [right] at (2.5,.7) {\scriptsize$\frac{\kappa-\rho}{1-\rho}$};
					\draw [dotted] (0,1.25) -- (1.25,1.25);
					\draw [thick] (-.07,1.25) -- (.07,1.25) node [left] at (0,1.25) {\scriptsize$\kappa$};
					\draw node [below] at (2,0) {\textcolor{white}{\scriptsize$\kappa/\rho$}};
					\draw node [above] at (1,2.5) {\textcolor{white}{\scriptsize$\frac{\kappa-(1-\rho)}{\rho}$}};	
				\end{tikzpicture}
			}
			&
			{\small$\kappa < \rho < 1$}
			\\\hline\\
		\end{tabular}
		\caption{Feasible sets of System \eqref{System:BasicModel} for $0 < \kappa < \frac{1}{2}$ and $\frac{1}{2} < \kappa < 1$.}
		\label{Figure:ActualFeasibleSets}
	\end{center}
\end{figure}
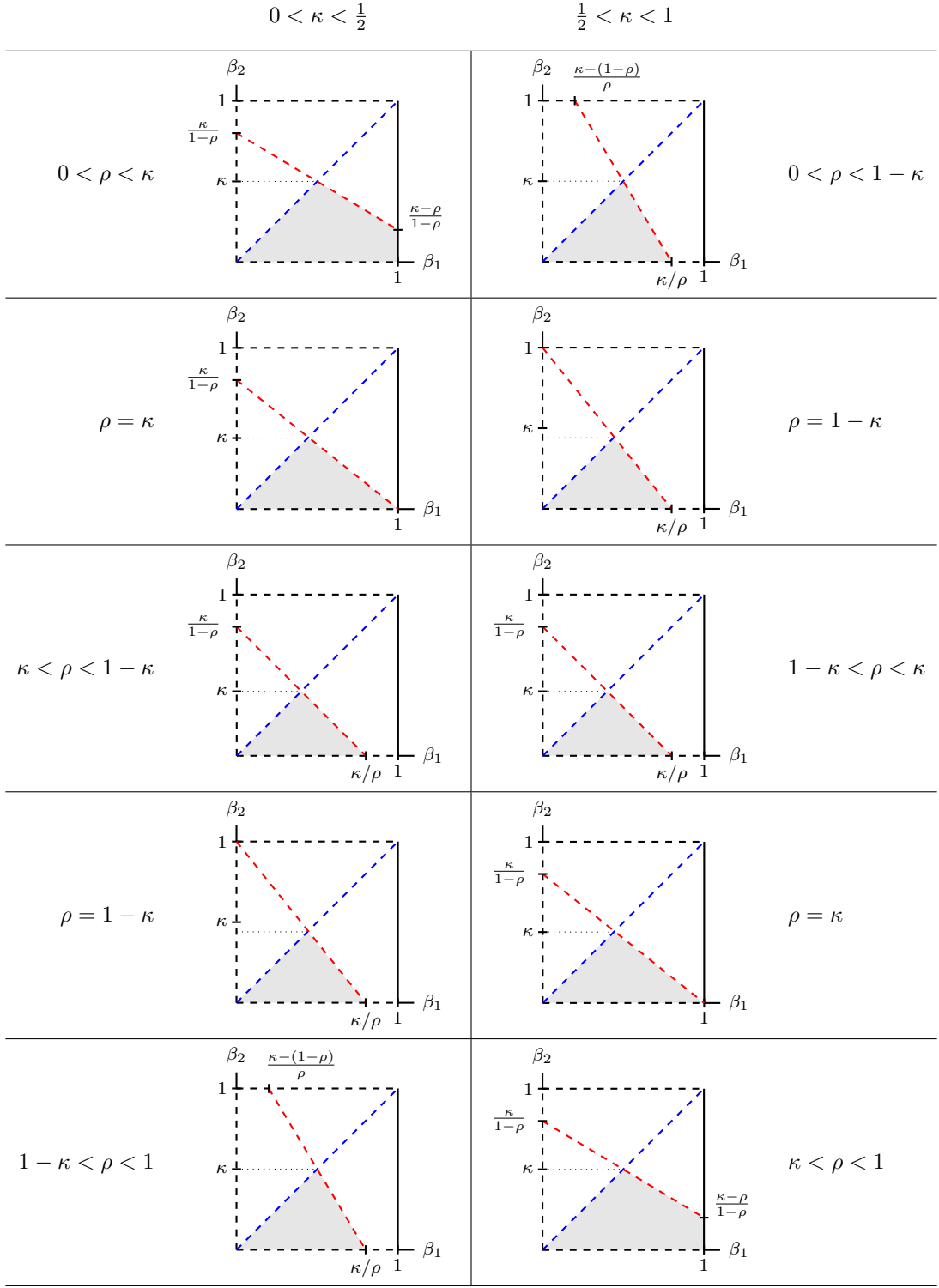

\begin{figure}[H]
	\begin{center}
		\begin{tabular}{c c c}
			\adjustbox{valign=c}{%
				\begin{tikzpicture}
					\fill [gray!20] (0,0) -- (1.25, 1.25) -- (2.5,.5) -- (2.5,0) -- (0,0) -- cycle;
					\draw [dashed, thick] (0,0) -- (0,2.5);
					\draw [thick] (0,2.5) -- (0,2.75);
					\draw node [above] at (0,2.75) {\scriptsize$\beta_2$};
					\draw [dashed, thick] (0,0) -- (2.5,0);
					\draw [thick] (2.5,0) -- (2.75,0);
					\draw node [right] at (2.75,0) {\scriptsize$\beta_1$};
					\draw [thick] (2.5,.07) -- (2.5,-.07) node [below] at (2.5,0) {\scriptsize$1$};
					\draw [thick] (.07,2.5) -- (-.07,2.5) node [left] at (0,2.5) {\scriptsize$1$};
					\draw [dashed, thick] (0,2.5) -- (2.5,2.5);
					\draw [thick] (2.5,0) -- (2.5,2.5);
					\draw [dashed, thick, blue] (0,0) -- (2.5,2.5);
					\draw [dashed, thick, red] (0,2) -- (2.5,.5);
					\draw [thick] (-.07,2) -- (.07,2) node [left] at (-.1,2) {\scriptsize$\frac{\kappa}{1-\rho}$};
					\draw [thick] (2.43,.5) -- (2.57,.5)  node [right] at (2.5,.7) {\scriptsize$\frac{\kappa-\rho}{1-\rho}$};
					\draw [dotted] (0,1.25) -- (1.25,1.25);
					\draw [thick] (-.07,1.25) -- (.07,1.25) node [left] at (0,1.25) {\scriptsize$\kappa$};
					\draw node [below] at (2,0) {\textcolor{white}{\scriptsize$\kappa/\rho$}};
					\draw node [above] at (.85,2.5) {\textcolor{white}{\scriptsize$\frac{\kappa-(1-\rho)}{\rho}$}};	
				\end{tikzpicture}   
			}
			&
			\adjustbox{valign=c}{%
				\begin{tikzpicture}
					\fill [gray!20] (0,0) -- (1.1, 1.1) -- (2.5,0) -- (0,0) -- cycle;
					\draw [dashed, thick] (0,0) -- (0,2.5);
					\draw [thick] (0,2.5) -- (0,2.75);
					\draw node [above] at (0,2.75) {\scriptsize$\beta_2$};
					\draw [dashed, thick] (0,0) -- (2.5,0);
					\draw [thick] (2.5,0) -- (2.75,0);
					\draw node [right] at (2.75,0) {\scriptsize$\beta_1$};
					\draw [thick] (2.5,.07) -- (2.5,-.07) node [below] at (2.5,0) {\scriptsize$1$};
					\draw [thick] (.07,2.5) -- (-.07,2.5) node [left] at (0,2.5) {\scriptsize$1$};
					\draw [dashed, thick] (0,2.5) -- (2.5,2.5);
					\draw [thick] (2.5,0) -- (2.5,2.5);
					\draw [dashed, thick, blue] (0,0) -- (2.5,2.5);
					\draw [dashed, thick, red] (0,2) -- (2.5,0);
					\draw [thick] (-.07,2) -- (.07,2) node [left] at (-.1,2) {\scriptsize$\frac{\kappa}{1-\rho}$};
					\draw [dotted] (0,1.1) -- (1.1,1.1);
					\draw [thick] (-.07,1.1) -- (.07,1.1) node [left] at (0,1.1) {\scriptsize$\kappa$};
					\draw node [below] at (2,0) {\textcolor{white}{\scriptsize$\kappa/\rho$}};
					\draw node [right] at (2.5,.7) {\textcolor{white}{\scriptsize$\frac{\kappa-\rho}{1-\rho}$}};
					\draw node [above] at (.85,2.5) {\textcolor{white}{\scriptsize$\frac{\kappa-(1-\rho)}{\rho}$}};	
				\end{tikzpicture}
			}
			& 
			\adjustbox{valign=c}{%
				\begin{tikzpicture}
					\fill [gray!20] (0,0) -- (1, 1) -- (2,0) -- (0,0) -- cycle;
					\draw [dashed, thick] (0,0) -- (0,2.5);
					\draw [thick] (0,2.5) -- (0,2.75);
					\draw node [above] at (0,2.75) {\scriptsize$\beta_2$};
					\draw [dashed, thick] (0,0) -- (2.5,0);
					\draw [thick] (2.5,0) -- (2.75,0);
					\draw node [right] at (2.75,0) {\scriptsize$\beta_1$};
					\draw [thick] (2.5,.07) -- (2.5,-.07) node [below] at (2.5,0) {\scriptsize$1$};
					\draw [thick] (.07,2.5) -- (-.07,2.5) node [left] at (0,2.5) {\scriptsize$1$};
					\draw [dashed, thick] (0,2.5) -- (2.5,2.5);
					\draw [thick] (2.5,0) -- (2.5,2.5);
					\draw [dashed, thick, blue] (0,0) -- (2.5,2.5);
					\draw [dashed, thick, red] (0,2) -- (2,0);
					\draw [thick] (-.07,2) -- (.07,2) node [left] at (-.1,2) {\scriptsize$\frac{\kappa}{1-\rho}$};
					\draw [dotted] (0,1) -- (1,1);
					\draw [thick] (-.07,1) -- (.07,1) node [left] at (0,1) {\scriptsize$\kappa$};
					\draw [thick] (2,-.07) -- (2,.07) node [below] at (2,0) {\scriptsize$\kappa/\rho$};
					\draw node [right] at (2.5,.7) {\textcolor{white}{\scriptsize$\frac{\kappa-\rho}{1-\rho}$}};
					\draw node [above] at (.85,2.5) {\textcolor{white}{\scriptsize$\frac{\kappa-(1-\rho)}{\rho}$}};	
				\end{tikzpicture}
			}\\
			{\small$0 < \rho < \kappa$} & {\small$\rho = \kappa$} & {\small$\kappa < \rho < \frac12$}\\
		\end{tabular}
		\caption{Feasible sets of System \eqref{System:ComplexModel} for $0 < \kappa < \frac12$.}\label{Figure:ComplexModel:ActualFeasibleSetsK<1/2}
	\end{center}
\end{figure}

\begin{figure}[H]
	\begin{center}
		\begin{tabular}{c c c}
			\adjustbox{valign=c}{%
				\begin{tikzpicture}
					\fill [gray!20] (0,0) -- (1.25, 1.25) -- (2,0) -- (0,0) -- cycle;
					\draw [dashed, thick] (0,0) -- (0,2.5);
					\draw [thick] (0,2.5) -- (0,2.75);
					\draw node [above] at (0,2.75) {\scriptsize$\beta_2$};
					\draw [dashed, thick] (0,0) -- (2.5,0);
					\draw [thick] (2.5,0) -- (2.75,0);
					\draw node [right] at (2.75,0) {\scriptsize$\beta_1$};
					\draw [thick] (2.5,.07) -- (2.5,-.07) node [below] at (2.5,0) {\scriptsize$1$};
					\draw [thick] (.07,2.5) -- (-.07,2.5) node [left] at (0,2.5) {\scriptsize$1$};
					\draw [dashed, thick] (0,2.5) -- (2.5,2.5);
					\draw [thick] (2.5,0) -- (2.5,2.5);
					\draw [dashed, thick, blue] (0,0) -- (2.5,2.5);
					\draw [dashed, thick, red] (.5,2.5) -- (2,0);
					\draw [thick] (.5,2.57) -- (.5,2.43) node [above] at (1,2.5) {\scriptsize$\frac{\kappa-(1-\rho)}{\rho}$};
					\draw [thick] (2,-.07) -- (2,.07) node [below] at (2,0) {\scriptsize$\kappa/\rho$};
					\draw [dotted] (0,1.25) -- (1.25,1.25);
					\draw [thick] (-.07,1.25) -- (.07,1.25) node [left] at (0,1.25) {\scriptsize$\kappa$};
					\draw node [left] at (-.1,2) {\textcolor{white}{\scriptsize$\frac{\kappa}{1-\rho}$}};
					\draw node [right] at (2.5,.7) {\textcolor{white}{\scriptsize$\frac{\kappa-\rho}{1-\rho}$}};
				\end{tikzpicture}
			}
			&
			\adjustbox{valign=c}{%
				\begin{tikzpicture}
					\fill [gray!20] (0,0) -- (1.1, 1.1) -- (2,0) -- (0,0) -- cycle;
					\draw [dashed, thick] (0,0) -- (0,2.5);
					\draw [thick] (0,2.5) -- (0,2.75);
					\draw node [above] at (0,2.75) {\scriptsize$\beta_2$};
					\draw [dashed, thick] (0,0) -- (2.5,0);
					\draw [thick] (2.5,0) -- (2.75,0);
					\draw node [right] at (2.75,0) {\scriptsize$\beta_1$};
					\draw [thick] (2.5,.07) -- (2.5,-.07) node [below] at (2.5,0) {\scriptsize$1$};
					\draw [thick] (.07,2.5) -- (-.07,2.5) node [left] at (0,2.5) {\scriptsize$1$};
					\draw [dashed, thick] (0,2.5) -- (2.5,2.5);
					\draw [thick] (2.5,0) -- (2.5,2.5);
					\draw [dashed, thick, blue] (0,0) -- (2.5,2.5);
					\draw [dashed, thick, red] (0,2.5) -- (2,0);
					\draw [thick] (2,-.07) -- (2,.07) node [below] at (2,0) {\scriptsize$\kappa/\rho$};
					\draw [dotted] (0,1.1) -- (1.1,1.1);
					\draw [thick] (-.07,1.25) -- (.07,1.25) node [left] at (0,1.25) {\scriptsize$\kappa$};
					\draw node [left] at (-.1,2) {\textcolor{white}{\scriptsize$\frac{\kappa}{1-\rho}$}};
					\draw node [above] at (1,2.5) {\textcolor{white}{\scriptsize$\frac{\kappa-(1-\rho)}{\rho}$}};
					\draw node [right] at (2.5,.7) {\textcolor{white}{\scriptsize$\frac{\kappa-\rho}{1-\rho}$}};
				\end{tikzpicture}
			}
			& 
			\adjustbox{valign=c}{%
				\begin{tikzpicture}
					\fill [gray!20] (0,0) -- (1, 1) -- (2,0) -- (0,0) -- cycle;
					\draw [dashed, thick] (0,0) -- (0,2.5);
					\draw [thick] (0,2.5) -- (0,2.75);
					\draw node [above] at (0,2.75) {\scriptsize$\beta_2$};
					\draw [dashed, thick] (0,0) -- (2.5,0);
					\draw [thick] (2.5,0) -- (2.75,0);
					\draw node [right] at (2.75,0) {\scriptsize$\beta_1$};
					\draw [thick] (2.5,.07) -- (2.5,-.07) node [below] at (2.5,0) {\scriptsize$1$};
					\draw [thick] (.07,2.5) -- (-.07,2.5) node [left] at (0,2.5) {\scriptsize$1$};
					\draw [dashed, thick] (0,2.5) -- (2.5,2.5);
					\draw [thick] (2.5,0) -- (2.5,2.5);
					\draw [dashed, thick, blue] (0,0) -- (2.5,2.5);
					\draw [dashed, thick, red] (0,2) -- (2,0);
					\draw [thick] (-.07,2) -- (.07,2) node [left] at (-.1,2) {\scriptsize$\frac{\kappa}{1-\rho}$};
					\draw [dotted] (0,1) -- (1,1);
					\draw [thick] (-.07,1) -- (.07,1) node [left] at (0,1) {\scriptsize$\kappa$};
					\draw [thick] (2,-.07) -- (2,.07) node [below] at (2,0) {\scriptsize$\kappa/\rho$};
					\draw node [right] at (2.5,.7) {\textcolor{white}{\scriptsize$\frac{\kappa-\rho}{1-\rho}$}};
					\draw node [above] at (1,2.5) {\textcolor{white}{\scriptsize$\frac{\kappa-(1-\rho)}{\rho}$}};	
				\end{tikzpicture}
			}
			\\
			{\small$0 < \rho < 1-\kappa$} & {\small$\rho = 1-\kappa$} & \small{$1-\kappa <\rho < \frac12$}
		\end{tabular}
		\caption{Feasible sets of System \eqref{System:ComplexModel} for $\frac12 <  \kappa < 1$.}\label{Figure:ComplexModel:ActualFeasibleSetsK>1/2}
	\end{center}
\end{figure}

\begin{figure}[H]
	\begin{center}
		\begin{tabular}{c c}
			\begin{tikzpicture}
				\draw [thick] (0,0) -- (3.25,0);
				\draw [thick] (0,-1) -- (0,1);
				\draw node [right] at (3.25,0) {\scriptsize$\rho$};
				\draw node [left, rotate=90] at (-.77,1.35) {\scriptsize Type of $\rho$-feasible set};
				\draw [thick] (1.5,.07) -- (1.5,-.07) node [below] at (1.5,0) {\scriptsize$\kappa$};
				\draw [thick] (3,.07) -- (3,-.07) node [below] at (3,0) {\tiny$\frac12$};
				\draw [thick] (.07,0) -- (-.07,0) node [left] at (0,0) {\scriptsize$0$};
				\draw [thick] (.07,.75) -- (-.07,.75) node [left] at (0,.75) {\scriptsize$1$};
				\draw [thick] (.07,-.75) -- (-.07,-.75) node [left] at (0,-.75) {\scriptsize$-1$};
				\draw [thick] (0,.75) -- (1.5,.75);
				\draw [thick] (1.5,-.75) -- (3,-.75);
				\filldraw [fill=white] (0,.75) circle (1.25pt);
				\filldraw [fill=white] (1.5,.75) circle (1.25pt);
				\filldraw [fill=white] (1.5,-.75) circle (1.25pt);
				\filldraw [fill=white] (3,-.75) circle (1.25pt);
				\filldraw (1.5,0) circle (1.25pt);
			\end{tikzpicture}
			&
			\begin{tikzpicture}
				\draw [thick] (0,0) -- (3.25,0);
				\draw [thick] (0,-1) -- (0,1);
				\draw node [right] at (3.25,0) {\tiny$\rho$};
				\draw node [left, rotate=90] at (-.77,1.35) {\scriptsize Type of $\rho$-feasible set};
				\draw [thick] (1.5,.07) -- (1.5,-.07) node [below] at (1.5,0) {\scriptsize$\kappa$};
				\draw [thick] (3,.07) -- (3,-.07) node [below] at (3,0) {\tiny$\frac12$};
				\draw [thick] (.07,0) -- (-.07,0) node [left] at (0,0) {\scriptsize$0$};
				\draw [thick] (.07,.75) -- (-.07,.75) node [left] at (0,.75) {\scriptsize$1$};
				\draw [thick] (.07,-.75) -- (-.07,-.75) node [left] at (0,-.75) {\scriptsize$-1$};
				\draw [thick] (0,-.75) -- (1.5,-.75);
				\draw [thick] (1.5,.75) -- (3,.75);
				\filldraw [fill=white] (0,-.75) circle (1.25pt);
				\filldraw [fill=white] (1.5,-.75) circle (1.25pt);
				\filldraw [fill=white] (1.5,.75) circle (1.25pt);
				\filldraw [fill=white] (3,.75) circle (1.25pt);
				\filldraw (1.5,0) circle (1.25pt);
			\end{tikzpicture}\\
			{\small$0<\kappa<\frac12$} & {\small$\frac12 < \kappa < 1$}
		\end{tabular}
	\end{center}
	\caption{Bifurcation diagrams of System \eqref{System:ComplexModel} as a function of $\rho$.}\label{Figure:BifurcationComplexModel}
\end{figure}

\newpage
\section*{Appendix B: Numerical Simulations of \eqref{System:BasicModel} and \eqref{System:ComplexModel}}

\newpage

The following simulations show the effects of varying the individual sociological parameters $\beta_1, \beta_2, \alpha_1, \alpha_2, \text{ and } \rho$.  Each of these simulations begin very close to an equilibrium.

\subsubsection*{Interesting Dynamics}
\begin{figure}[H]
\begin{picbox}
	\includegraphics[width=\textwidth]{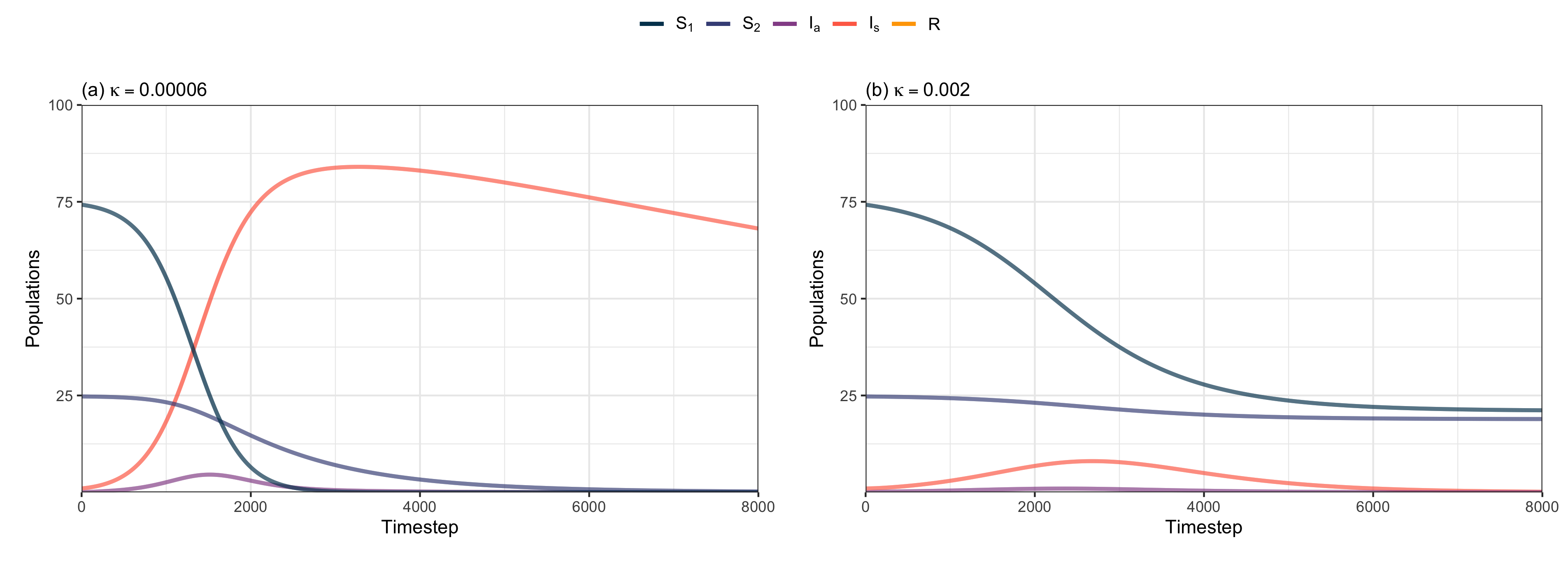}
	\tcblower
	\begin{tabular}{x{.035\textwidth}x{.15\textwidth}x{.15\textwidth}x{.15\textwidth}x{.15\textwidth}x{.15\textwidth}}
		{\footnotesize (a)} & $\beta_1$={\footnotesize 0.0042} & $\beta_2$={\footnotesize 0.0009} &  $\gamma$={\footnotesize 0.005} & $\lambda$={\footnotesize 0.65} & $\kappa$={\footnotesize 0.00006}
		\\
		{\footnotesize (b)} & $\beta_1$={\footnotesize 0.0042} & $\beta_2$={\footnotesize 0.0009} &  $\gamma$={\footnotesize 0.005} & $\lambda$={\footnotesize 0.65} & $\kappa$={\footnotesize 0.00200}
    \end{tabular}
\end{picbox}
\caption{Simulations of System \eqref{System:BasicModel} with $\rho = 0.75$ and initial condition $(S_1 = 99\rho, S_2 = 99(1-\rho), I_s = 1, I_a = 0, R = 0)$.}\label{Figure:BasicModelSimulations}
\end{figure}

\begin{figure}[H]
\begin{picbox}
	\includegraphics[width=\textwidth]{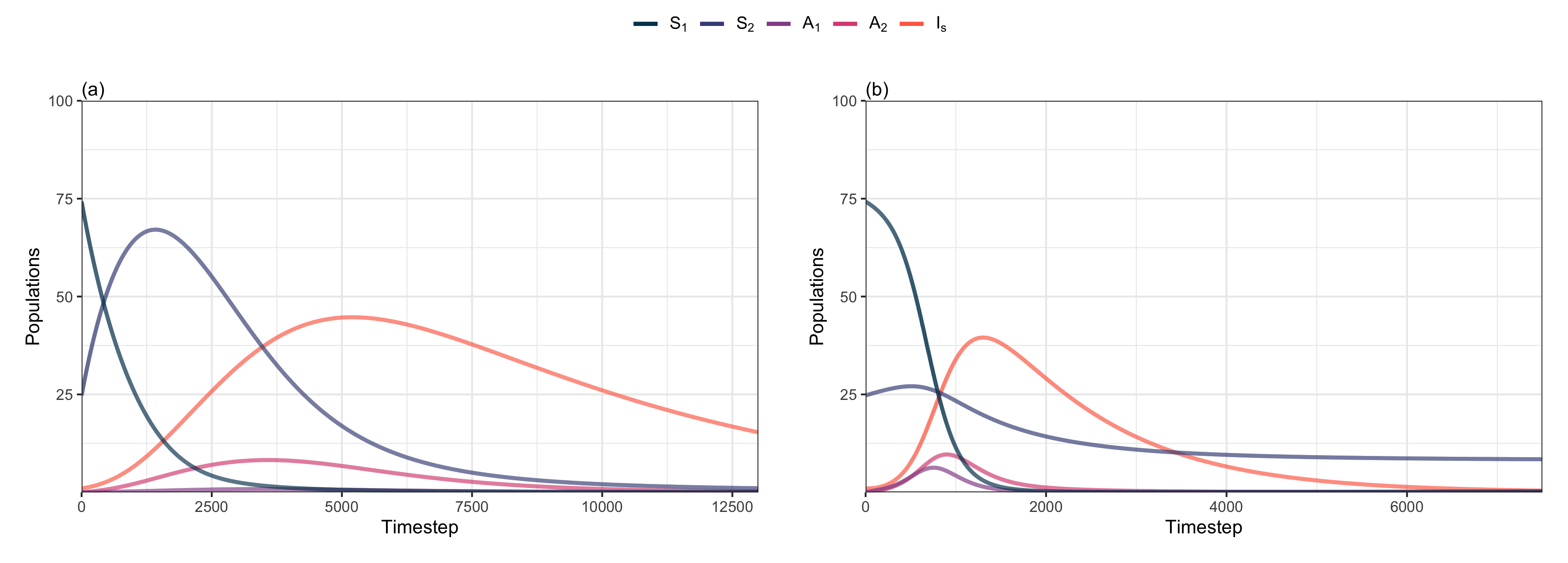}
	\tcblower
	\noindent\begin{tabular}{lx{.12\textwidth}x{.12\textwidth}x{.1\textwidth}x{.1\textwidth}x{.1\textwidth}x{.1\textwidth}x{.1\textwidth}}
	    (a) & $\beta_1$={\footnotesize 0.0042} & $\beta_2$={\footnotesize 0.0009} & $\alpha_1$={\footnotesize 0.10} & $\alpha_2$={\footnotesize 0.010} & $\gamma$={\footnotesize 0.0005} & $\lambda$={\footnotesize 0.65} & $\kappa$={\footnotesize 0.0002}
	    \\
        (b) & $\beta_1$={\footnotesize 0.0090} & $\beta_2$={\footnotesize 0.0012} & $\alpha_1$={\footnotesize 0.01} & $\alpha_2$={\footnotesize 0.001} & $\gamma$={\footnotesize 0.0050} & $\lambda$={\footnotesize 0.05} & $\kappa$={\footnotesize 0.0009}
	\end{tabular}
\end{picbox}	
\caption{Simulations of System \eqref{System:ComplexModel} with $\rho = 0.75$ and initial condition $(S_1 = 99\rho, S_2 = 99(1-\rho), I_s = 1, I_a = 0, R = 0)$.}\label{Figure:ComplexModelSimulations}
\end{figure}

Figures \ref{Figure:BasicModelSimulations} and \ref{Figure:ComplexModelSimulations} demonstrate that both models can exhibit extreme end behavior.

\newpage

\subsubsection*{Changing parameter $\beta_1$ and $\beta_2$}
\begin{figure}[H]
\begin{picbox}
	\includegraphics[width=\textwidth]{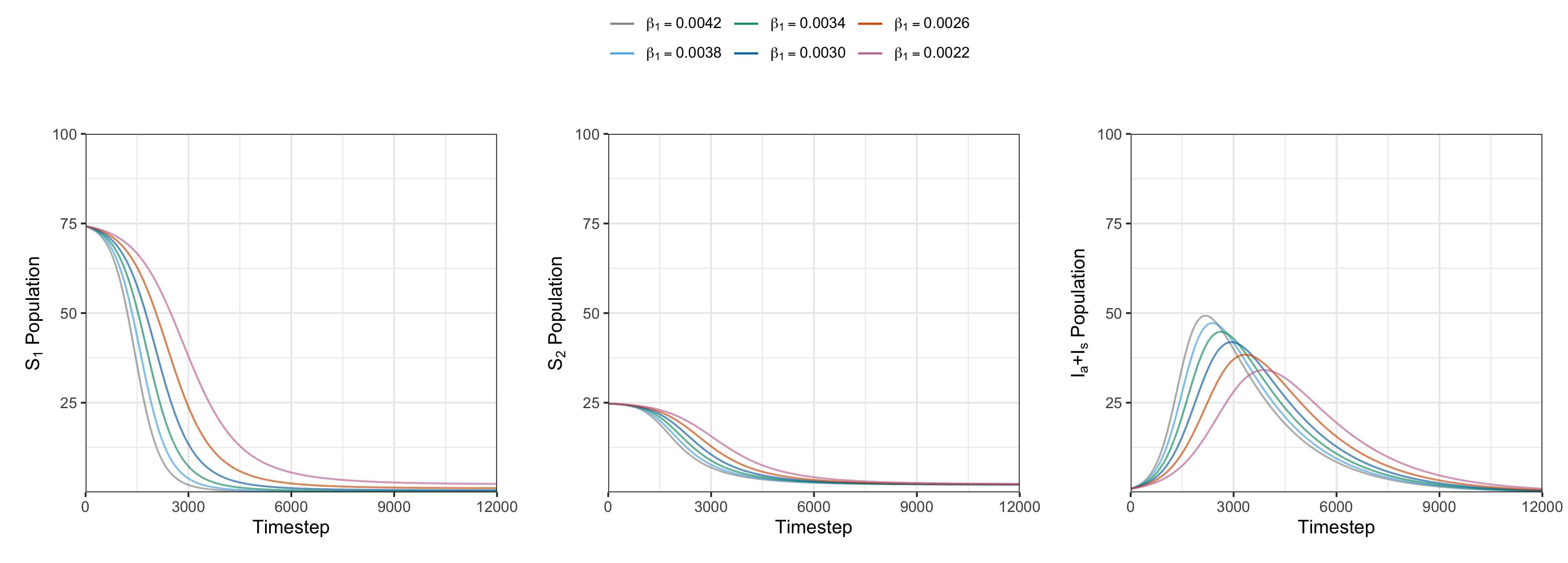}
	\tcblower
	\begin{tabular}{x{.17\textwidth}x{.17\textwidth}x{.17\textwidth}x{.17\textwidth}x{.17\textwidth}}
		$\beta_2$ = {\footnotesize 0.0015} & $\kappa$ = {\footnotesize 0.0006} & $\gamma$ = {\footnotesize 0.005} & $\lambda$ = {\footnotesize 0.65} & $\rho$ = {\footnotesize 0.75}
	\end{tabular}
\end{picbox}
\caption{Populations $S_1, S_2,$ and $I_a+I_s$ in Model \eqref{System:BasicModel} as $\beta_1$ varies from $0.0042$ to $0.0022$ with initial condition $(S_1 = 99\rho, S_2 = 99(1-\rho), I_s = 1, I_a = 0, R = 0)$.}\label{Figure:BasicModelEffectOfB1}
\end{figure}

\begin{figure}[H]
\begin{picbox}
	\includegraphics[width=\textwidth]{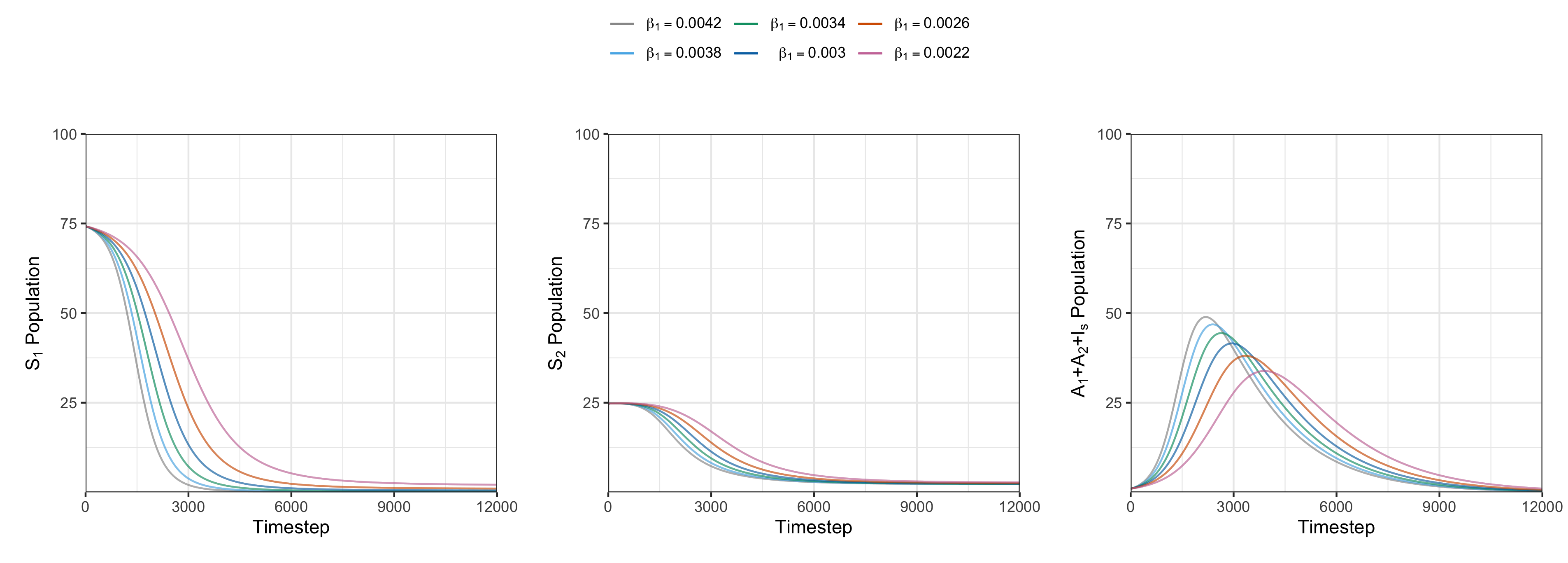}
	\tcblower
	\begin{tabular}{x{.115\textwidth}x{.115\textwidth}x{.115\textwidth}x{.1\textwidth}x{.1\textwidth}x{.115\textwidth}x{.12\textwidth}}
		$\beta_2$ = {\scriptsize 0.0015} & $\kappa$ = {\footnotesize 0.0006} & $\gamma$ = {\footnotesize 0.005} & $\lambda$ = {\footnotesize 0.05} & $\rho$ = {\footnotesize 0.75} & $\alpha_1$ = {\footnotesize 0.001} & $\alpha_2$ = {\scriptsize 0.0001}
	\end{tabular}
\end{picbox}
\caption{Populations $S_1, S_2,$ and $A_1+A_2+I_s$ in Model \eqref{System:ComplexModel} as $\beta_1$ varies from $0.0042$ to $0.0022$ with initial condition $(S_1 = 99\rho, S_2 = 99(1-\rho), A_1 = 0, A_2 = 0, I_s = 1, R = 0)$.}\label{Figure:ComplexModelEffectOfB1}
\end{figure}

Figures \ref{Figure:BasicModelEffectOfB1} and \ref{Figure:ComplexModelEffectOfB1} show the effect of changing the parameter $\beta_1$.  

\newpage

\begin{figure}[H]
\begin{picbox}
	\includegraphics[width=\textwidth]{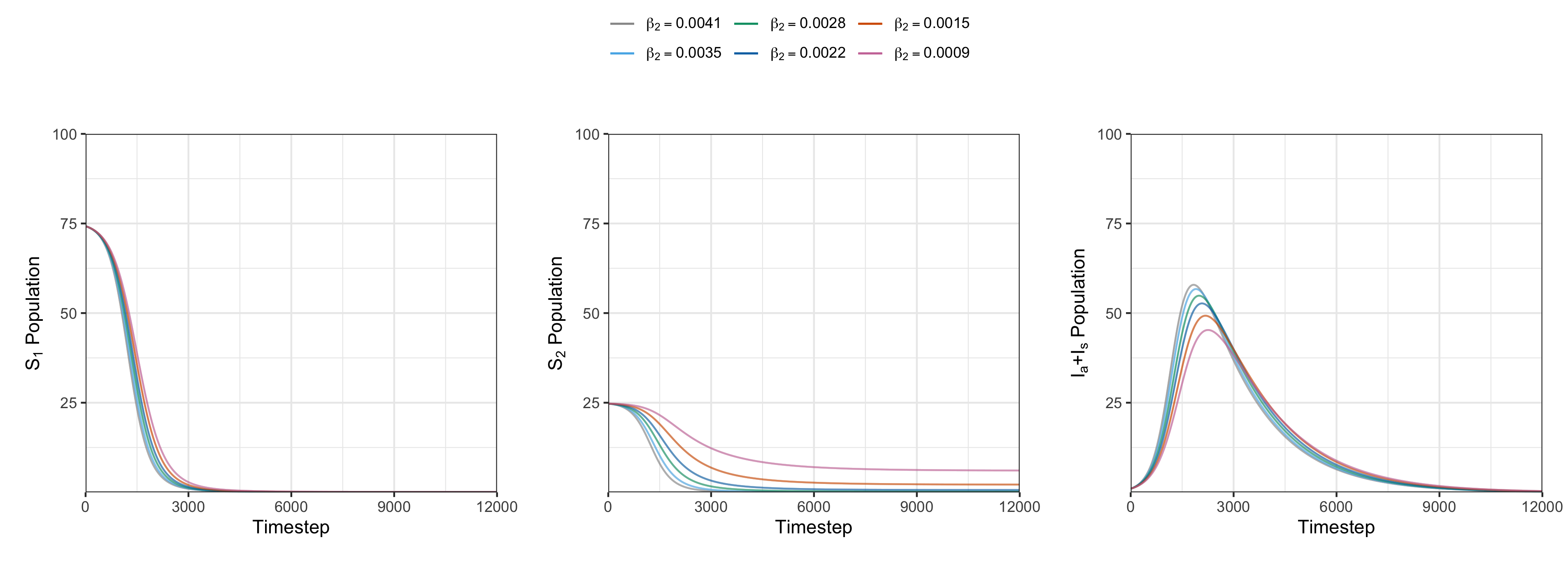}
	\tcblower
	\begin{tabular}{x{.17\textwidth}x{.17\textwidth}x{.17\textwidth}x{.17\textwidth}x{.17\textwidth}}
	    $\beta_1$ = {\footnotesize 0.0042} & $\kappa$ = {\footnotesize 0.0006} & $\gamma$ = {\footnotesize 0.005} & $\lambda$ = {\footnotesize 0.65} & $\rho$ = {\footnotesize 0.75}
	\end{tabular}
\end{picbox}
\caption{Populations $S_1, S_2,$ and $I_a+I_s$ of Model \eqref{System:BasicModel} as $\beta_2$ varies from $0.0041$ to $0.0009$ with initial condition $(S_1 = 99\rho, S_2 = 99(1-\rho), I_s = 1, I_a = 0, R = 0)$.}\label{Figure:BasicModelEffectOfB2}
\end{figure}

\begin{figure}[H]
\begin{picbox}
	\includegraphics[width=\textwidth]{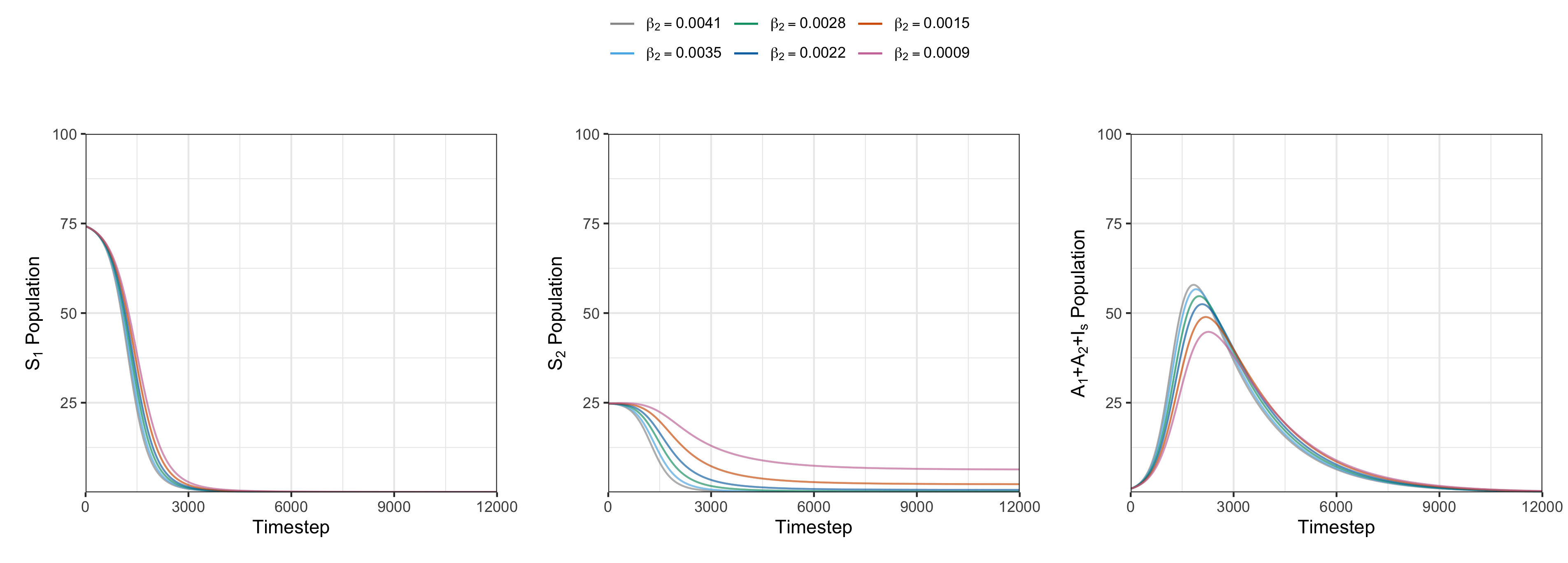}
	\tcblower
	\begin{tabular}{x{.115\textwidth}x{.115\textwidth}x{.115\textwidth}x{.1\textwidth}x{.1\textwidth}x{.115\textwidth}x{.12\textwidth}}
	    $\beta_1$ = {\scriptsize 0.0042} & $\kappa$ = {\footnotesize 0.0006} & $\gamma$ = {\footnotesize 0.005} & $\lambda$ = {\footnotesize 0.65} & $\rho$ = {\footnotesize 0.75} & $\alpha_1$ = {\footnotesize 0.001} & $\alpha_2$ = {\scriptsize 0.0001}
	\end{tabular}
\end{picbox}
\caption{Populations $S_1, S_2,$ and $A_1+A_2+I_s$ of Model \eqref{System:ComplexModel} as $\beta_2$ varies from $0.0041$ to $0.0009$ with initial condition $(S_1 = 99\rho, S_2 = 99(1-\rho), A_1 = 0, A_2 = 0, I_s = 1, R = 0)$.}\label{Figure:ComplexModelEffectOfB2}
\end{figure}

\noindent Figures \ref{Figure:BasicModelEffectOfB2} and \ref{Figure:ComplexModelEffectOfB2} show the effect of changing the parameter $\beta_2$.

\noindent It is worth noting that a casual inspection of the above figures might suggest that both systems may have a similar sensitivity to changes in $\beta_1$ and $\beta_2$.  This notion is supported by the sensitivity analysis in Section 3.3 and the discussion within.

\newpage
\subsubsection*{Changing parameters $\alpha_1$ and $\alpha_2$}

\begin{figure}[H]
\begin{picbox}
	\includegraphics[width=\textwidth]{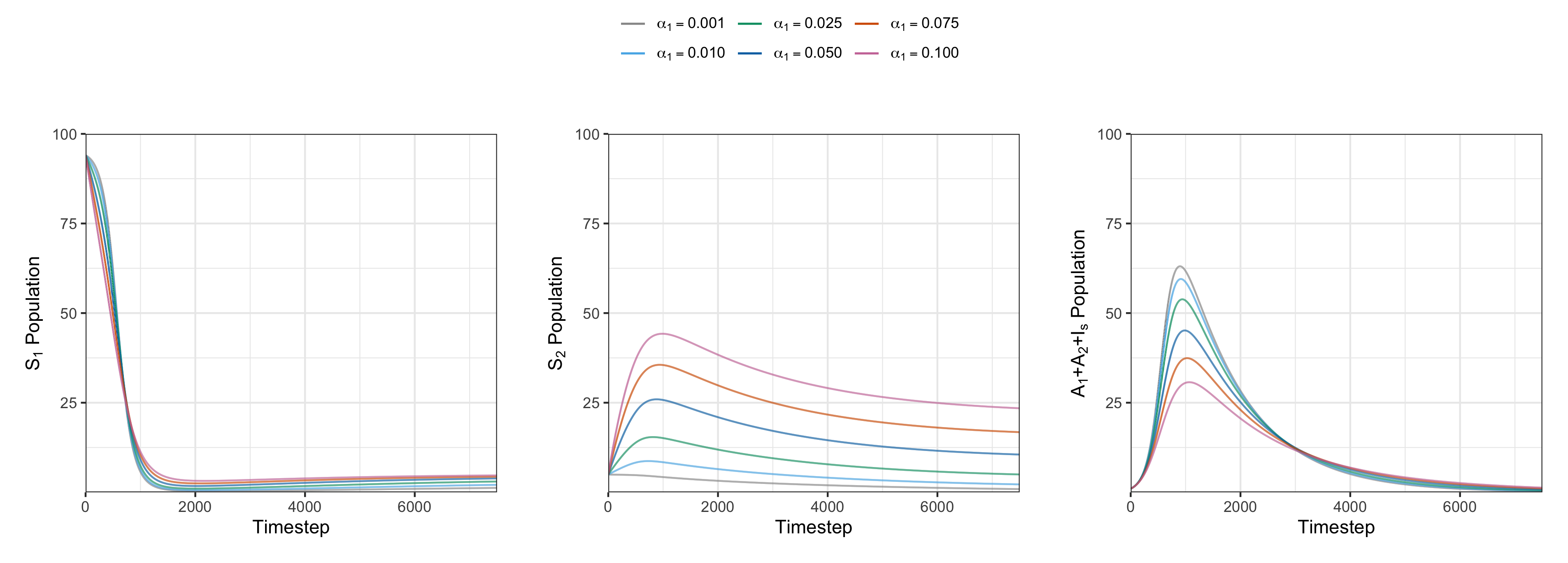}
    \tcblower
	\begin{tabular}{x{.135\textwidth}x{.135\textwidth}x{.135\textwidth}x{.135\textwidth}x{.135\textwidth}x{.135\textwidth}}
	    $\beta_1$ = {\footnotesize 0.009} & $\beta_2$={\footnotesize 0.00012} & $\kappa$ = {\footnotesize 0.0009} & $\gamma$ = {\footnotesize 0.005} & $\lambda$ = {\footnotesize 0.05} & $\alpha_2$ = {\footnotesize 0.023}
	\end{tabular}	
\end{picbox}
\caption{Populations $S_1, S_2,$ and $A_1+A_2+I_s$ of Model \eqref{System:ComplexModel} as $\alpha_1$ varies from $0.001$ to $0.1$ with initial condition $(S_1 = 99\rho, S_2 = 99(1-\rho), A_1 = 0, A_2 = 0, I_s = 1, R = 0)$.}\label{Figure:ComplexModelEffectOfA1}
\end{figure}

Figure \ref{Figure:ComplexModelEffectOfA1} shows the impact of changing the parameter $\alpha_1$ on System \eqref{System:ComplexModel}.

\begin{figure}[H]
\begin{picbox}
	\includegraphics[width=\textwidth]{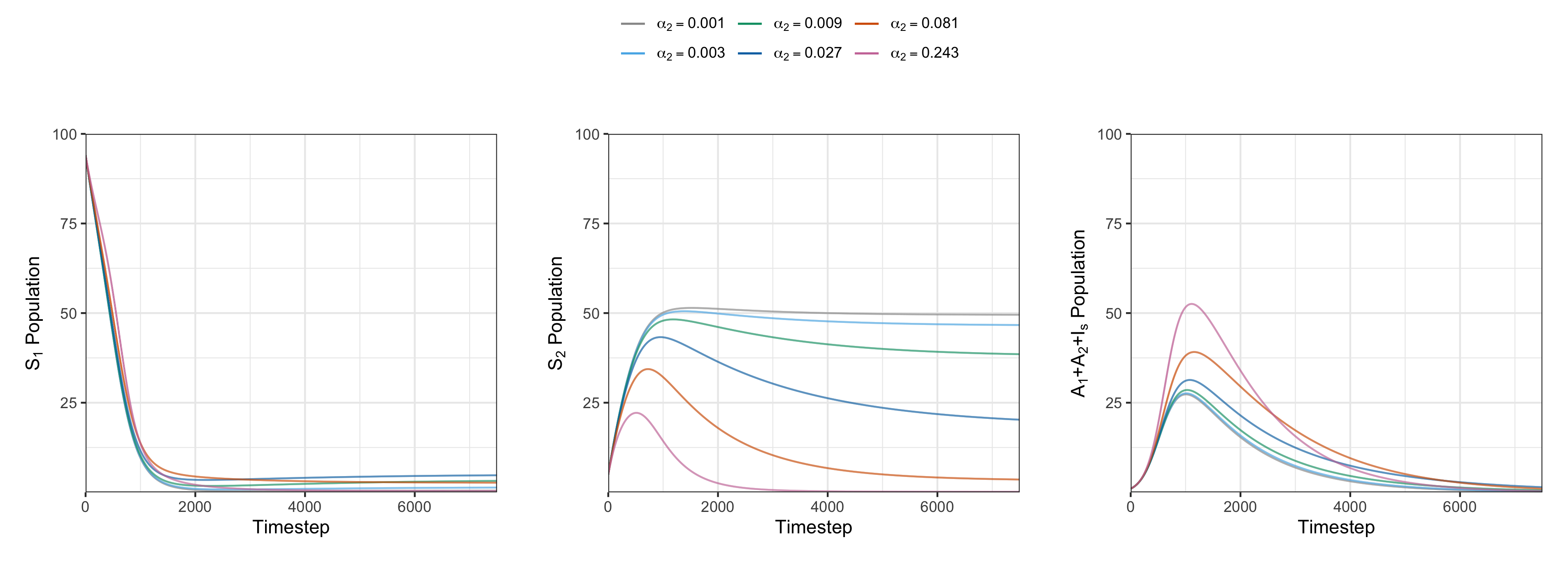}
    \tcblower
	\begin{tabular}{x{.135\textwidth}x{.135\textwidth}x{.1435\textwidth}x{.135\textwidth}x{.135\textwidth}x{.135\textwidth}}
	    $\beta_1$ = {\scriptsize 0.009} & $\beta_2$={\scriptsize 0.00012} & $\kappa$ = {\scriptsize 0.0009} & $\gamma$ = {\footnotesize 0.005} & $\lambda$ = {\footnotesize 0.05} & $\alpha_1$ = {\scriptsize 0.1} 
	\end{tabular}	
\end{picbox}
\caption{Populations $S_1, S_2,$ and $A_1+A_2+I_s$ of Model \eqref{System:ComplexModel} as $\alpha_2$ varies from $0.001$ to $0.243$ with initial condition $(S_1 = 99\rho, S_2 = 99(1-\rho), A_1 = 0, A_2 = 0, I_s = 1, R = 0)$.}\label{Figure:ComplexModelEffectOfA2}
\end{figure}

Figure \ref{Figure:ComplexModelEffectOfA2} shows the impact of changing the parameters $\alpha_1$ on System \eqref{System:ComplexModel}. It is worth noting here that smaller values of $\alpha_2$ not only result in a smaller peak in infections, but that peak occurs slightly earlier compared with larger values of $\alpha_2$.

\newpage
\subsubsection*{Changing parameter $\rho$}
\begin{figure}[H]
\begin{picbox}
    \includegraphics[width=\textwidth]{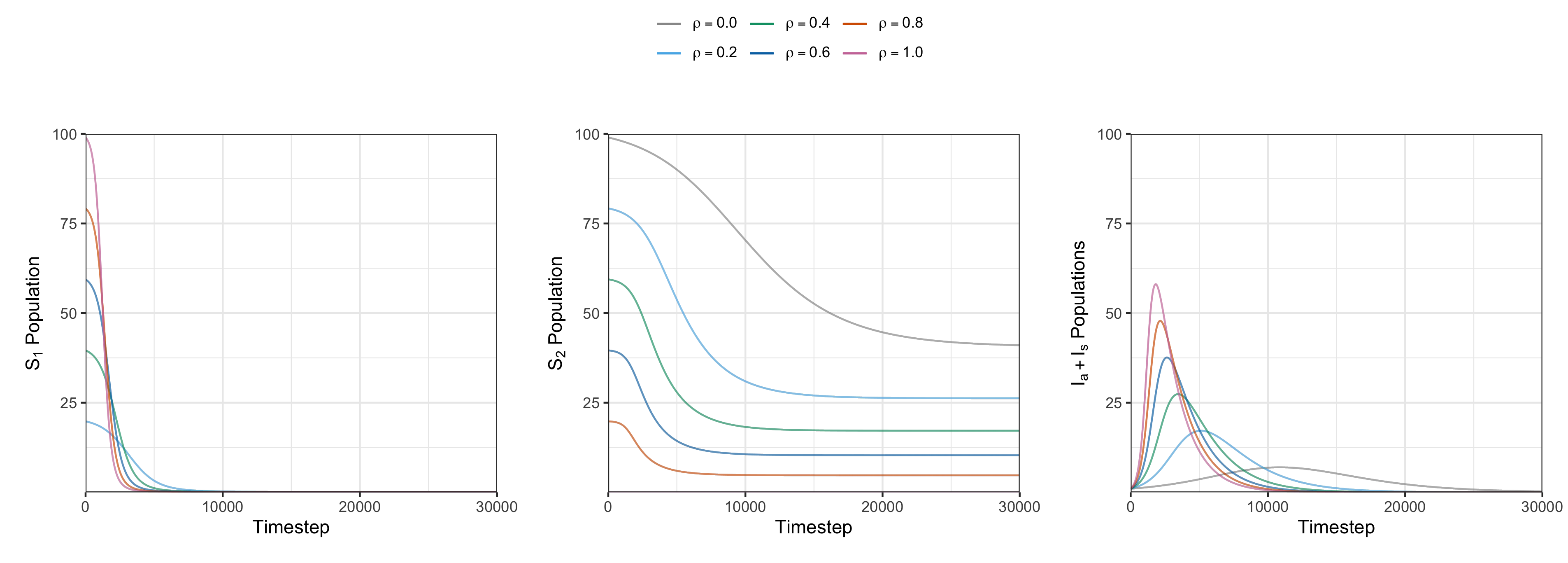}
	\tcblower
	\begin{tabular}{x{.17\textwidth}x{.17\textwidth}x{.17\textwidth}x{.17\textwidth}x{.17\textwidth}}
        $\beta_1$ = {\footnotesize 0.0042} & $\beta_2$ = {\footnotesize 0.0009} & $\kappa$ = {\footnotesize 0.0006} & $\gamma$ = {\footnotesize 0.005} & $\lambda$ = {\footnotesize 0.65}
	\end{tabular}
\end{picbox}
\caption{Populations $S_1, S_2,$ and $I_a+I_s$ in Model \eqref{System:BasicModel} as $\rho$ varies from $0$ to $1$ with initial condition $(S_1 = 99\rho, S_2 = 99(1-\rho), I_s = 1, I_a = 0, R = 0)$.}\label{Figure:BasicModelChangeinRho}
\end{figure}

Figure \ref{Figure:BasicModelChangeinRho} shows the effect of changing the parameter $\rho$ on System \eqref{System:BasicModel}.  There is no such figure for System \eqref{System:ComplexModel} as $\rho$ is not a free parameter.

\end{document}